\documentclass[12pt]{amsart}
\usepackage{times}      
\usepackage{latexsym}   
\usepackage{amssymb}    
\usepackage{amsmath}    
\usepackage{amsbsy}
\usepackage{amsthm}
\usepackage{amsgen}
\usepackage{amsfonts}
\usepackage{array}
\usepackage{epsfig}

\newtheorem{theorem}{Theorem}[section]
\newtheorem*{thm0}{Theorem}
\newtheorem*{thma}{Theorem A}
\newtheorem*{thmb}{Theorem B}
\newtheorem*{thmeq}{Equivariant Sphere Theorem}
 
\newtheorem{conjecture}[theorem]{Conjecture} 
\newtheorem*{conjc}{Conjecture C}
\newtheorem{lemma}[theorem]{Lemma} 
\newtheorem{proposition}[theorem]{Proposition}

\def\rrr{\mathbb{R}}
\def\ccc{\mathbb{C}}
\def\zzz{\mathbb{Z}}

\def\hh{\mathbb{H}}

\DeclareMathOperator{\diam}{diam}
\DeclareMathOperator{\lcm}{lcm}
\DeclareMathOperator{\Out}{Out}
\DeclareMathOperator{\Isom}{Isom}
\DeclareMathOperator{\Fix}{Fix}
\def\bdm{\begin{displaymath}}
\def\edm{\end{displaymath}}
\def\beq{\begin{equation}}
\def\eeq{\end{equation}}
\def\bes{\begin{equation*}}
\def\ees{\end{equation*}}
\def\epcm{\end{picture}\end{center}\end{minipage}}
\def\bpcm{\begin{minipage}{80pt}\begin{center}\begin{picture}}
\def\so{SO}

\def\f4{F_4}
\def\g2{G_2}
\def\normal{\vartriangleleft}
\def\mo{\mathbf{O}}
\def\mt{\mathbf{T}}
\def\mi{\mathbf{I}}
\def\mii{\mathbf{i}}
\def\mj{\mathbf{j}}
\def\mk{\mathbf{k}}
\def\md{\mathbf{D}}
\def\mc{\mathbf{C}}

\begin{document} 

\title[Diameters of 3-Sphere Quotients]{Diameters of 3-Sphere Quotients}

\author[Dunbar]{William D. Dunbar}
\address[Dunbar]{Department of Mathematics, 
Simon's Rock College, Great Barrington, MA  01230, U. S. A.}
\email{wdunbar@simons-rock.edu}

\author[Greenwald]{Sarah J. Greenwald$^*$}
\address[Greenwald]{Department of  Mathematics, Appalachian 
State University, Boone, NC 28608, U. S. A.}
\email{greenwaldsj@appstate.edu}

\author[McGowan]{Jill McGowan}
\address[McGowan]{Department of Mathematics, Howard University,
Washington, DC 20059, U. S. A.}
\email{jmcgowan@fac.howard.edu}

\author[Searle]{Catherine Searle$^{**}$}
\address[Searle]{Institute of Mathematics, University Nacional Autonoma 
de Mexico, Cuernavaca, Morelos, MEXICO}
\email{csearle@matcuer.unam.mx}


\thanks{$^*$
The second author was supported in part by NSF ROA grants 0072533 and 9972304.
$^{**}$The fourth author was supported in part by 
CONACYT Project \#SEP-CO1-46274.}


\begin{abstract}

Let $G\subset O(4)$ act isometrically on $S^3$. In this article we calculate
a lower bound for the diameter of the quotient spaces $S^3/G$. We
find it to be 
$\frac{1}{2}\arccos(\frac{\tan(\frac{3 \pi}{10})}{\sqrt3})$,
which is exactly the value of the lower bound for diameters
of the spherical space forms. In the process, we are also able to find a 
lower bound for diameters for the spherical Aleksandrov spaces, $S^n/G$,
of cohomogeneities 1 and 2, as well as for cohomogeneity 3 (with some 
restrictions on the group type). This leads us to conjecture that
the diameter of $S^n/G$ is increasing as the cohomogeneity of 
the group $G$ increases.
\end{abstract}
\maketitle
\section{Introduction}\label{S:intro}
Diameter is one of the most basic geometric invariants.  Knowing its 
lower bound not only provides information about the orbit space $X=S^n/G$, but 
also leads to other interesting results. 
And, while representations of compact Lie groups are well understood, the geometry 
of the corresponding spherical quotients is virtually unknown and is 
potentially very important.  Let  $X^k=S^n/G$, where $G$ is a closed, 
nontransitive subgroup of $O(n)$, and examine the diameter of $X^k$.  
When G is finite, $k=n$ and $X^n$ is a manifold, there exits an explicit lower bound 
on the diameter that
depends only on the dimension $n$, and a global lower bound  that is 
independent of dimension also exists \cite{M93}. (Throughout this paper, unless otherwise explicitly stated, $S^n$ is taken to be the round unit sphere of dimension $n$.) For other closed, 
nontransitive groups, a lower bound on the diameter also exists \cite{Gr00}, 
but it is not always given explicitly.  This paper examines the diameter of 
quotients of the three-dimensional sphere to find an optimal lower bound. 
We also find descriptions of many of these orbit spaces 
and learn about their geometry.

Often, given a specific compact 
Lie group, $G$, acting isometrically and (almost) effectively on $M$, with 
$sec(M)>0$, we can recover the manifold from the possible orbit space 
decompositions of the 
action. For example, a classic theorem of
Hsiang and Kleiner \cite{HK89} states:
\begin{thm0}\label{T:HsiangKleiner}
Let $M^4$ be a 1-connected, strictly 
positively curved closed Riemannian manifold which admits an effective,
isometric $S^1$ action. Then $M^4$ is homeomorphic to $S^4$ or $CP^2$.
\end{thm0}
Here, the essential point of the proof lies in understanding
the fixed point set of the circle action. By work of Freedman \cite{F82}, one can 
recover the manifold merely by knowing 
the Euler characteristic,  $\chi(M)$,  
which in this case turns out to be equal 
to  $\chi(\Fix(M;S^1))$. In particular, in the case where
$\Fix(M;S^1)$ consists of isolated points, these points are
singular orbits of the action, and we can bound the total number of such 
orbits via the diameter of the orbit space of the normal sphere to 
any such point of isotropy. Here the normal sphere is an $S^3$
and the upper bound on the diameter of $S^3/S^1$ tells us that there 
are no more than 3 such points.
We note as well, that a theorem of Rong \cite{R02}
uses the same technique to show that a 1-connected, strictly positively curved,
closed Riemannian
5-manifold admitting a $T^2$ isometric and effective action
is homeomorphic to $S^5$. 
In particular, this is part of a more general phenomenon
where a bound on the {\it q-extent} of a space allows us to limit
the number of singular points of 
a given action, and we have (\cite{GM95}, \cite{GS97}):
\begin{thmeq}
Let $M$ be a closed manifold with $sec(M)>0$ on which $G$ acts 
(almost) effectively by isometries. 
Suppose $p_0,p_1\in M$ are points such that 
$\diam S_{\bar{p_i}}\leq \pi/4$, $i= 0,1$, 
where $S_{\bar{p_i}}$ is the space of directions 
at $\bar{p_i}$ in $M/G$. Then $M$ can be exhibited as 
\[
M=D(G(p_0))\bigcup_E D(G(p_1))
\]
where $D(G(p_i))$, $i=0,1$ are tubular neighborhoods of the 
$p_i$-orbits and $E=\partial D(G(p_0))=\partial D(G(p_1))$. In 
particular, $M$ is homeomorphic to the sphere if $G(p_i)=p_i$, i.e., if 
$p_i$, $i=0,1$ are isolated fixed points of $G$ and  
$\diam S_{\bar{p_i}}\leq \pi/4$.
\end{thmeq}
Thus, local diameter information gives 
global results about the structure of the manifold.  

\medskip

When $G$ is finite in $O(n+1)$, then $S^n/G$ is a good orbifold, that is, 
the global quotient of a Riemannian manifold by  a discrete subgroup of its 
isometry group \cite{B91}.  
Finite subgroups of $O(4)$ are classified in \cite{D64} and various methods 
from McGowan \cite{M93} and Dunbar \cite{Du94} are used to find a lower bound 
on the diameter of the resulting spherical quotients.  When $G$ is infinite, 
$S^n/G$ is a spherical Alexandrov space with curvature bounded below.  This is 
a length space with Riemannian notions such as distance and curvature obtained 
by comparison with $S^n$ via Toponogov \cite{BGP92}.  Possible groups in 
$SO(4)$ are classified in 
\cite{S96} and \cite{S97}.  Extensions of these groups in $O(4)$ are examined 
along with the 
diameters of the resulting spherical quotients.
\begin{thma}
If $G$ is a closed, non-transitive subgroup in $O(4)$ then
\[
\diam(S^3 /  G ) \geq  \frac{\alpha}{2}
\]
where $\alpha=\arccos(\frac{\tan(\frac{3 \pi}{10})}{\sqrt3})$.
\end{thma}
This diameter is approximately 
$\frac{\pi}{9.63}$ and is achieved by $S^3 / \eta(S^1 \times \mathbf{I})$, where 
$\mathbf{I}$ is the binary icosahedral group, and 
$\eta: Sp(1) \times Sp(1) \rightarrow SO(4)$ is defined by first noting that the unit quaternions may be identified with $S^3$ by 
$\phi( p_1 + ip_2 + jp_3 + kp_4)=(p_1,p_2,p_3,p_4)$, where $p_1^2 + p_2^2 + p_3^2 + p_4^2= 1$. 
With this identification, $\eta$ maps $(a,b)$ to $A$, if for every $x \in Sp(1)$, 
$\phi(axb^{-1}) = A\phi(x)$.  The map $\eta$ is a surjective homomorphism with kernel
$\{(1,1),(-1,-1)\}$.
 
If $G$ is finite then
$S^3 / \eta(\mathbf{C}_{2m} \times \mathbf{I})$, in the limit as $m\rightarrow \infty$, achieves 
the smallest diameter, 
where $\mathbf{C}_{2m}$ is the binary cyclic group.   The orbit space is a 
manifold only if $\gcd(m, 30)=1$ \cite{W84} and 
otherwise it is a Seifert-fibered orbifold, foliated by
circles and intervals.  Among nonfibering orbifolds,  
$S^3/\eta(\mathbf{O} \times \mathbf{I})$ achieves the 
smallest diameter, where $\mathbf{O}$ is 
the binary octahedral group.  A lower bound estimate for this diameter is 
$ \arccos\left( 1/(\sqrt{40+12\sqrt{2}-8\sqrt{5}-12\sqrt{10}})\right)$, which is 
approximately $\frac{\pi}{8.93}$.

Note that  the {\it cohomogeneity} of a connected 
$G$-action is the codimension of its principal
orbit, or equivalently, the dimension of the orbit space. We 
extend this definition to include non-trivial disconnected actions.
With respect to cohomogeneity one actions, we not only examine those actions
on $S^3$, but also on a round sphere of any dimension, and we find
that the smallest diameter for a non-trivial disconnected cohomogeneity one
action on any $S^n$ is $\pi/6$.
We show that only certain actions of 
cohomogeneity one admit a finite extension of the group 
that halves the diameter of the corresponding orbit space.

Using this result and results from \cite{MS05} and 
\cite{MS06} we are able to prove the 
following:
\begin{thmb}
 Let $G$ act by cohomogeneity 1, 2, or 3 on $S^n$. Further, suppose 
that if $G$ is connected, the action is a classical connected polar action. 
Then
\bdm
\min(\diam(S^n/G))=\left\{ \begin{array}{ll}
\frac{\pi}{6} & \textrm{for cohomogeneity 1}\\
\frac{\alpha}{2} & \textrm{for cohomogeneity 2}\\
\frac{\alpha}{2} & \textrm{for cohomogeneity 3}
\end{array} \right.
\edm

\end{thmb}
Here we define a {\it classical polar}
action to be one which corresponds to a symmetric space $G/H$
where both $G$ and $H$ are classical Lie groups. Recall that 
all polar actions on spheres correspond to the isotropy 
subgroup of a symmetric space \cite{D85}, i.e., they correspond to
the natural action of the isotropy subgroup $H$ of $G/H$ acting on 
$T_{G(e)}G/H$.

This theorem and other work (\cite{MS05}, \cite{MS06}) 
lead us to the following conjecture:
\begin{conjc}
Let $G$ act irreducibly on $S^n$ by cohomogeneity $k$, where $n\in 2\zzz$.  Then 
for all $\epsilon > 0$ and for all sufficiently large $k$
(and all $n > k$), $\diam(S^n/G)$ is within $\epsilon$ of $\pi/2$.
\end{conjc}
 We break the remainder of the paper into three sections, which each
consider the action of $G$ on $S^3$ of a specific cohomogeneity, plus
a section giving our conclusions.

\section{Cohomogeneity One}\label{S:cohomOne}
\begin{proposition}
Let $G$ be an action on $S^3$ whose 
orbit space is dimension 1, 
then the minimal diameter of $S^3/G$ is $\pi/4$.
\end{proposition}
\begin{proof}
The classification of low cohomogeneity actions on spheres 
(\cite{S96, S97}, \cite{HL71}) 
tells us that the only two possible connected groups which can act 
effectively on $S^3$ by 
cohomogeneity one are $SO(3)$ and $T^2$.

The action by  $SO(3)$  has principal orbit 
$SO(3)/SO(2)\simeq S^2$, and singular orbits equal to points 
($SO(3)/SO(3)$). Its orbit space is an interval of length $\pi$.
The second action is $T^2$ acting with principal orbit $T^2$ and singular 
orbits $T^1$ (each singular orbit is a different $T^1$). Its orbit space is 
an interval of length $\pi/2$. 

\medskip

\noindent {\bf Observation}:
{\it Let $H'$ be a finite extension of $H$, a connected subgroup of $G$, 
then $H'\subset N_G(H)$.}

\smallskip

Note that conjugation is a group isomorphism and that 
given any $g\in H'$, $gHg^{-1}$ is isomorphic to $H$, the connected component
of the identity of $H'$. Further, since $gHg^{-1}$ is also a subgroup of $H'$,
$gHg^{-1}=H$, and $H'\subset N_G(H)$.

\medskip

 We may write the action of $SO(3)$ acting 
by cohomogeneity one on $S^3$, 
as follows. Let $S^3$ be the standard sphere in 
$\rrr^4$; we represent its points as $(x,y,z,w)$ and the action is
represented by the matrix 
\[
A=\left( \begin{matrix}
B & 0\\
0 & 1\\
\end{matrix} \right),
\]
where $B\in SO(3)$.
It is clear that the points $(0,0,0,1)$ and $(0,0,0,-1)$ are fixed by this 
action.
As well, it is clear that the matrix 
\[
C=\left( \begin{matrix}
1 & 0 & 0 & 0\\
0 & 1 & 0 & 0\\
0 & 0 & 1 & 0\\
0 & 0 & 0 & -1\\
\end{matrix} \right)
\]
is an element of $O(4)$, and that  $C$ commutes with $A$ and acts on $S^3$ 
by interchanging
the $w$ coordinate of any point with its negative.
One easily sees that this action sends any given orbit with a specified $w$ 
coordinate to the corresponding orbit with $-w$ coordinate.
This action ``folds'' the orbit space of $S^3/SO(3)$ 
(an interval of length $\pi$) in half to obtain an interval of length $\pi/2$.
Moreover, the normalizer of $SO(3)$ in $O(4)$ is exactly the
group generated by $SO(3)$ and $C$.

Now $T^2$ is the maximal torus 
in $O(4)$ and we will denote it by
\[
T^2=\left(\begin{matrix}
\cos(\theta_1) & \sin(\theta_1) & 0 & 0\\
-\sin(\theta_1) & \cos(\theta_1) &0 &0 \\
0 & 0 & \cos(\theta_2) & \sin(\theta_2)\\
0 & 0 & -\sin(\theta_2) & \cos(\theta_2)\\
\end{matrix}\right).
\]
The Weyl group of $SO(4)$ is isomorphic to $C_2\times C_2$, where
$C_n$ is the cyclic group of order $n$,
and the order 2 element
\[
\left(
\begin{matrix}
0 & 0 & 1 & 0\\
0 & 0 & 0 & 1\\
1 & 0 & 0 & 0\\
0 & 1 & 0 & 0\\
\end{matrix}
\right)
\]
identifies the two singular orbits
$T^2/S^1$ to each other, and acts upon the principal orbit $T^2$
at distance $\pi/4$ from the singular orbits, identifying
the remaining principal orbits in pairs according to their
respective distances from the singular orbits.
Thus, the diameter of the resulting orbit space is $\pi/4$.

We claim that any other finite extension cannot further reduce the 
diameter of the orbit space due to the following lemma,
which will complete the proof.
\end{proof}
\begin{lemma}\label{T:cuthalf}
Suppose $G$ acts isometrically and (almost)
effectively on $M^n$ by cohomogeneity one,
where $M^n$ is a closed Riemannian manifold of strictly positive 
sectional curvature. Then any finite extension of $G$ in $\Isom(M^n)$
can reduce the diameter of $M^n/G$ by at most one-half.
\end{lemma}
\begin{proof}
Topologically, there are 4 possibilities
for the orbit space of a manifold by a cohomogeneity one action. They are
$\rrr$, $\rrr^+$, $S^1$, and an interval (cf. \cite{M57}). 
However, the additional restriction of strictly positive sectional curvature
eliminates the first three possibilities, since such a manifold will be compact
and have finite fundamental group.
Over the interval, the manifold decomposes into principal orbits over its
interior and 2 singular orbits over each endpoint, and the diameter
of the orbit space is given by the length of the interval.

Now, any finite isometric action on such an interval can only ``fold'' 
the interval in half, identifying the endpoints to each other 
and corresponding pairs of principal orbits at a given distance from
one of the corresponding pairs of endpoints.
The principal orbit equidistant from both endpoints is not identified to
any other, and itself is acted on nontrivially. 
Clearly such an action halves the diameter.
Note that any other action would be discontinuous or not isometric.
Since any finite group of order greater than 2 will have a cyclic element
of order greater than or equal to 3, or a subgroup isomorphic to $C_2\times C_2,$ it suffices 
to understand these two cases.  

In the case where we have a subgroup of order 3 acting effectively,
the action would identify three distinct points in the interval. However,
it would fail 
to be isometric 
(suppose $a\neq b$, $b\neq c$, $c\neq a$, and $f(a)=b$, 
$f(b)=c$, and $f(c)=a$, then we must 
have $d(a,b) = d(f(a),f(b)) = d(f(b),f(c))$
in order for it to be an isometry, but this is clearly impossible in any 
interval).  The argument is similar if we have a cyclic group of order 
greater than three.

For the case where we have an effective action by a subgroup isomorphic to $C_2\times C_2$, 
if the action decreases the diameter, then it must fold the interval in 
half two times. For the action to be isometric, this means that 
the corresponding singular orbits of the original interval and 
of the once-folded interval must respectively be isometric themselves.
We observe however that even if this is true at first, it cannot hold for the 
once-folded interval, since one singular orbit of this interval will
be one of the original singular orbits, and the other 
the result of a $\mc_2$ action on a principal orbit.
A quick inspection of Table \ref{Ta:cohomOne} of cohomogeneity one actions
on spheres shows that in all cases, the principal orbit is of
strictly larger dimension than its 
singular orbits and thus 
it is impossible to make any further
identifications, that is, an effective action by a finite subgroup isomorphic to $C_2\times C_2$
is not allowed.
\end{proof}
\begin{proposition}\label{T:CohomOneBound}
The smallest diameter 
one can obtain for the orbit space of a (non-trivial, disconnected) 
cohomogeneity one 
action on a sphere is $\pi/6$.
\end{proposition}
\begin{proof}
We will show that 
the cohomogeneity one actions of diameter $\pi$ and $\pi/2$
both admit finite extensions which fold the corresponding interval in half.
As well, there are two actions of diameter $\pi/3$ (of the total four)
which also admit finite extensions. 
The $\pi$ and $\pi/2$ diameter actions 
are exactly the reducible actions.
The remaining irreducible actions, other 
than the two of diameter $\pi/3$ we mentioned previously, do not admit
finite extensions. We
note that all the orbit spaces of diameter $\pi/4$ have non-isometric 
singular orbits and thus there is no isometric action which can reduce their
respective diameters (see proof of Lemma \ref{T:cuthalf}). Thus we need only show
that the remaining $\pi/3$ and the $\pi/6$ actions do not 
admit finite extensions
which halve the diameters of their orbit spaces.

 In Table \ref{Ta:cohomOne} (cf. \cite{MS05}, \cite{S97}),
we give a list of the cohomogeneity one actions on spheres. 
\begin{table}
  \begin{center}
  \caption{Connected Spherical Cohomogeneity One Actions}\label{Ta:cohomOne}
    \begin{tabular}{lccc}
 Group$(G)$ & Representation($\Phi)$ & $dim(\Phi)$ & Length   \\ 
1) $SO(k)$ &$\rho_k+1$, $k\geq 2$ & $k+1$& $\pi$ \\
2) $U(k)$ & $\mu_k+1$, $k\geq 1$ & $2k+1$ & $\pi$\\
3) $Sp(k)$ & $\nu_k+1$, $k\geq 1$ & $4k+1$ & $\pi$ \\
4) $G_2$ & $\psi_1 +1$ & $7$ & $\pi$ \\
5) $Spin(7)$ & $\Delta_7 +1$ & $8$& $\pi$\\
6) $Spin(9)$& $\Delta_9 +1$ &$16$ & $\pi$\\
7) $SO(k) \times SO(m)$& $\rho_k+\rho_m$, $k,m\geq 2$ & $k+m$ & $\pi/2$ \\
8) $SO(3)$ & $ S^2\rho_3-1$ & $5$ & $\pi/3$  \\
9) $SU(3)$&  Ad &$8$&  $\pi/3$\\
10) $Sp(3)$ &$\wedge^2\nu_3-1$ & $14$& $\pi/3$\\ 
11) $F_4$ & $\phi_1$ & $26$ & $\pi/3$\\ 
12) $SO(2) \times SO(k)$ &$\rho_2 \otimes_{\rrr} \rho_k$, 
$k\geq 3$ &$2k$ & $\pi/4$\\
13)  $U(2) \times SU(k)$ &  $\mu_2 \otimes_{\ccc} \mu_k$, 
$k\geq 2$ & $4n$&$\pi/4$\\
14) $Sp(2) \times Sp(k)$ & $\nu_2 \otimes_{\hh}\nu_k$, 
$k\geq 2$ & $8n$ & $\pi/4$\\
15) $U(5)$& $[\wedge^2\mu_5]_{\rrr}$ & $20$ & $\pi/4$\\
16) $Sp(2)$ & Ad & $10$&$\pi/4$\\
17) $U(1) \times Spin(10)$ &  
$[\mu_1 \otimes_{\ccc}\Delta^{\pm}_{10}]_{\rrr}$ &$32$ & $\pi/4$\\
18) $G_2$&  Ad & $14$&  $\pi/6$ \\
19)  $SO(4)$ & $\nu_1 \otimes_{\hh} S^3 \nu_1$ & $8$&  $\pi/6$ \\
    \end{tabular}
  \end{center}
\end{table}
Note that in order for an isometric action to actually fold 
the orbit space in half, the two singular orbits 
(as mentioned in the proof of Lemma \ref{T:cuthalf})
must be isometric themselves, since they will be identified to 
each other via an isometry. 
In the list of cohomogeneity one spherical orbit spaces,  
only the following have isometric singular orbits: numbers
1, 2, 3, 4, 5, 6, 7 (with $k=m$), 8, 9, 10, 11, 17, 18, and 19.
The first three
 all have diameters which reduce to $\pi/2$
with the addition of an antipodal action.
For example, in case 1, where $SO(k)$ acts on a sphere of dimension 
$k$, we see that adding in the element 
\[
\left(\begin{matrix}
I_{k\times k} & 0\\
0 & -1\\
\end{matrix}
\right)
\]
which belongs to $O(k+1)=\Isom(S^{k})$ identifies the two singular orbits
(which are points here) to each other, 
each of 
the two principal orbits equidistant from each of the singular orbits
will identify and the action on the principal orbit halfway between
the
two singular orbits is antipodal. Note that the action here is reducible 
and not maximal. Clearly this element also conjugates 
the two singular isotropy subgroups (which in this case
are the entire group $SO(k)$) to one another.
The remaining cases 2--6  with diameter $\pi$ proceed in a similar fashion.

In number 7, we may add in the element
\[
\left(\begin{matrix}
0 & I_{k\times k}\\
I_{k\times k} & 0\\
\end{matrix}\right)
\]
which is an element of $O(2k)=\Isom(S^{2k-1})$. This action
interchanges the two singular orbits (here they are $S^{k-1}$'s)
and corresponding principal orbits. The principal orbit equidistant
from the 2 singular orbits is acted upon antipodally and we obtain 
an orbit space of diameter $\pi/4$.
Note that this action is maximal, but it is reducible.
Further note that both singular istropy subgroups are
conjugate to each other precisely via 
this order 2 element.

Numbers 8, 9, 10, 11 all have diameter $\pi/3$ and the principal orbits
are flag manifolds and the corresponding singular orbits
are projective spaces (respectively real, complex, hyperbolic and Cayley).
Numbers 18 and 19 both have diameter $\pi/6$.
Of these actions, only numbers 8 and 9 admit a finite extension which
will fold the interval in half.
Note that all these actions correspond to maximal inclusions of
their respective groups in the corresponding
isometry group of the sphere on which they act by cohomogeneity one.

In the particular case of number 8, 
the action may be described as follows: 
$\so(3)$ acts on $S^4$ realized as 
symmetric $3\times 3$ 
real matrices of trace zero by conjugation. If we require that the symmetric 
matrix $A$ also satisfies $\|A\|^2={\rm tr} (A^tA)=1,$ then this is an action 
on a subset of $S^8$. These matrices can be diagonalized by the action, 
thus the orbits can be represented 
by the diagonal matrices with the appropriate eigenvalues whose sum is zero.
Further, conjugation by the matrices 
\bdm
\left(\begin{matrix}
0 & -1& 0\\
1& 0 & 0\\
0& 0& 1\\
\end{matrix}\right)
\textrm{and} 
\left(\begin{matrix}
1 & 0& 0\\
0& 0 & -1\\
0& 1& 0\\
\end{matrix}\right),
\edm 
 allows us to arrange the eigenvalues 
in descending order:  $x \geq y \geq z$. 
The resulting orbit space is 
then the subset of the intersection of $S^2$ with the plane 
$\{(x,y,z) \in S^2 : x+y+z=0\}$ with $x \geq y \geq z$, i.e., the segment of 
the great circle in $S^2$ with endpoints 
$({\frac{1}{\sqrt6}},{\frac{1}{\sqrt6}},{\frac{-2} {\sqrt6}})$ and 
$({\frac{2}{\sqrt6}},{\frac{-1} {\sqrt6}},{\frac{-1} {\sqrt6}})$.  No 
further identifications can be made by conjugation, since conjugation of 
matrices does not change their eigenvalues. 

However, multiplying on the left and on the right by the matrix 
\[
A=
\left(\begin{matrix}
i & 0& 0\\
0& i & 0\\
0& 0& i\\
\end{matrix}\right),
\]
gives us an antipodal map, which interchanges the endpoints of the interval and effectively halves the diameter 
of the resulting orbit space.

The action of $SU(3)$ is similar and we may likewise conjugate by the matrix
\[
A=
\left(\begin{matrix}
j & 0& 0\\
0& j & 0\\
0& 0& j\\
\end{matrix}\right).
\]
Note that this map is not the antipodal map on the corresponding sphere, in fact it is not even fixed point free.
However, it does map orbits to orbits (one can see this via the corresponding 
eigenvalues of a given orbit)
and acts antipodally on the geodesic circle containing the orbit space, and thus on the orbit space 
itself.  Therefore this action also will halve the orbit space diameter.

Now, for the remaining cases 10, 11,
18 and 19, we note that all of the actions $G$ correspond to 
a maximal inclusion of $G$ in the corresponding isometry group of 
the sphere upon which they act. That is 
$Sp(3)\subset O(14), F_4\subset O(28),\so(4)\subset O(8)$ and
$G_2\subset O(14)$.
Thus, if there exists such an element $g$ which finitely extends
the group $G$, $g\in O(k)\diagdown G$, then
we claim that  $g\in Aut(G)$. We see this as follows: we know that $g(G)$ must also be a subgroup of
$O(k)$ since it must take isotropy subgroups of 
$G$ to each other. Further, since 
$G$ is maximally included, then  $g(G)=G$, for if
they were not equal, then by maximality we know that 
their union would be all of $O(k)$, which clearly 
contradicts the assumption that 
$g$ was an element that extends the group finitely,
and 
Thus, $g\in Aut(G)$, that is, $g$ must be contained in the automorphism group of the action
$G$. 
In particular, the inner automorphisms
will not give us extensions (finite or otherwise), whereas the outer ones might.
Thus those groups which have trivial outer automorphism group will
not admit finite extensions halving the diameter. 
In the remaining cases, it is well-known that only
for $SO(4)$ the outer automorphism group is non-trivial.
We will now show that $SO(4)$ does not admit such a finite extension.

From Uchida \cite{U80}, we have a description of the action on $S^7\subset \rrr^8\simeq \hh^2.$  We begin with the homomorphism $\sigma: SU(2)\rightarrow Sp(2)$
defined by 
\[
\sigma\left(\begin{matrix}
a & -\bar{b}\\
b & \bar{a}\\
\end{matrix}\right)= \left(\begin{matrix}
a^3 + jb^3 & -\sqrt{3}(a^2\bar{b} - j\bar{a}b^2)\\
\sqrt{3}(a^2b - jab^2) & a^2\bar{a}-2ab\bar{b} + jb^2\bar{b} -2ja\bar{a}b\\
\end{matrix}\right),
\]
where $j$ is a quaternion such that $j^2=-1$ and $aj=j\bar{a}$ for each 
complex number $a$.

Next, we let $A\in \sigma(SU(2)),\ q\in Sp(1),\ X\in M(2, 1;\hh)$,
then the action is defined by $(A, q)\times X\mapsto AX \bar q$.

Now $\Out(SO(4))\cong \mathbf C_2$ and can be generated via the action
$\tau: Sp(1)\times Sp(1)\rightarrow Sp(1)\times Sp(1)$, where
$\tau(q,r)=(r,q)$.
Since Uchida describes the action using the double cover 
$Sp(1)\times Sp(1)$ of $SO(4),$ we will work with the double cover as well. (See the description  of the double cover in the introduction.)
We will suppose that an extension by $\tau$ exists and derive a contradiction.
Suppose there exists $\beta\in Sp(1)\times Sp(1)$
such that $\beta\tau$ interchanges orbits as desired (it suffices to consider
this case since $\tau\beta\tau^{-1}\in Sp(1)\times Sp(1)$ implies 
that $\tau\beta=\gamma\tau$ for some $\gamma\in Sp(1)\times Sp(1)$.  Now
given
that $\beta\tau$ is such an extension, in particular
it maps elements of $G(0, j)$ to elements of $G(\frac{1}{\sqrt{3}}, j)$.
Observe that 
\bdm
G_{(0,j)}=(\sigma\biggl(\left(\begin{matrix}
a & 0\\
0 & \bar{a}\\
\end{matrix}\right)\biggr), \bar{a}),
\edm
and that $\tau G_{(0,j)}\tau^{-1}=G_{(0,j)}$

By a direct calculation, one can see that the orbit of 
$(\frac{1}{\sqrt{3}}, j)$ contains no elements of the form
$(0, bj),\, b \in \ccc$, which are the only elements in 
$S^7$ that are fixed by $G_{(0,j)}$. Thus $\tau$ cannot send the 
point $(0, j)$ to any point of the other singular orbit 
$G(\frac{1}{\sqrt{3}}, j)$ (if it did then the isotropy subgroup 
would be $G_{(0,j)}$, since its image under conjugation by $\tau$ is itself.)
Thus $\tau(0,j)\in G_{(0,j)}$.
Likewise $\beta\tau(0,j)\in G_{(0,j)}$. Thus there exists no finite 
extension that folds the interval in half.

\end{proof}
%
\section{Cohomogeneity Two}\label{S:cohomTwo}

\begin{figure}
\begin{center}
\epsfxsize=70mm \epsfbox{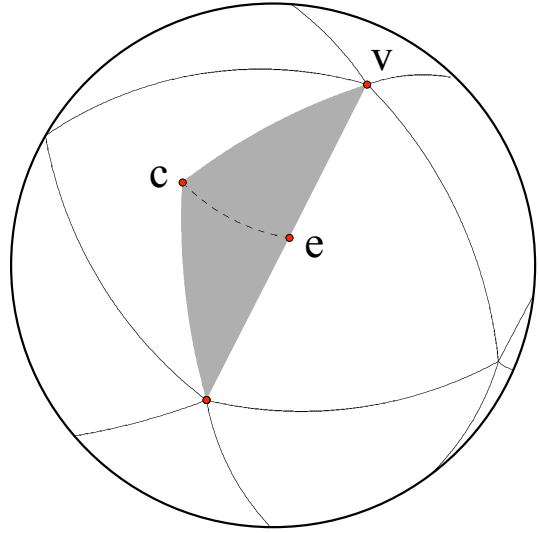}
\end{center}
\caption{Fundamental Domain of $I$ on the Spherical Icosahedron}
\label{Fig:233}
\end{figure}

In this section we compute a lower bound on the diameters of $S^3/G$ of 
dimension two.
We begin by proving two lemmas.
\begin{lemma}\label{T:componentG}
If $G$ is a topological group acting 
by isometries on the $n$-sphere $S^n$,
and $G_0$ denotes the connected component of the identity
element, then $G/G_0$ acts on  $S^n/G_0$ by 
isometries.
\end{lemma}
\begin{proof}
First, note that $S^n/G_0$ has a metric 
that is well 
defined and the distance between two orbits is given by  
\begin{align}
d(G_0x,G_0y) &= \min_{g,h \in G_0} d(gx,hy)=d(g_0x,h_0y), 
{\rm for \,\, some \,}g_0, \, h_0 \in G_0 \notag \\
& = d(h_0^{-1}g_0x,y) \geq d(G_0x,y) \geq d(G_0x,G_0y) \notag
\end{align}
Now, let $G/G_0$ act on $S^n/G_0$.  Let $hG_0 \in G/G_0$.  We want to 
determine $d(hG_0(G_0x),hG_0(G_0y))$.  But this is equal to 
\begin{align}
d(hG_0x,hG_0y)&=d(G_0hx,G_0hy)=d(G_0G_0hx,hy)=d(G_0hx,hy) \notag\\
&=d(h^{-1}G_0hx,y)=d(G_0x,y)=d(G_0x,G_0y) \notag
\end{align}  
So the action by $G/G_0$ preserves distance.
\end{proof}
Hence, in order to examine the diameter of $S^n/G$, we will look at the 
diameter of
$(S^n/G_0)/ (G/G_0)$, which is isometric to $S^n/G$  \cite{M93}.
In fact, Lemma \ref{T:closedO3} proves more than we actually need 
for Proposition \ref{T:CohomTwoBound}, since it considers 
the question of any closed subgroup of $O(3)$ acting on $S^2$.
\begin{lemma}\label{T:closedO3}
If $G$ is a closed, non-transitive subgroup in $O(3)$ then 
\[
\diam(S^2 / G ) \geq  \alpha,
\]
where $\alpha=\arccos(\frac{\tan(\frac{3 \pi}{10})}{\sqrt3})$.
\end{lemma}
\begin{proof}
If $G$ acts by cohomogeneity one, from Table \ref{Ta:cohomOne} we see that the 
only possible connected action on $S^2$ 
is by $\so(2)$ and gives diameter $\pi$. From
the proof of Proposition \ref{T:CohomOneBound} we see that the smallest diameter 
we can obtain with a finite extension of $\so(2)$ is $\frac{\pi}{2}$
and this occurs for $G=O(2)O(1)$.

If $G$ acts by cohomogeneity 2, then $G$ is finite, and 
the possibilities for $G$ and the 
corresponding diameters are listed in the following 
Table \ref{Ta:2Dfin}. For each group, we also provide 
several standard notations.  
The smallest diameter is  achieved by $S^2/ I$ and $S^2/ I^h$
as $\arccos{\frac{\tan(\frac{3 \pi}{10})}{\sqrt3}}$ \cite{Gr98}, where 
$I$ is the icosahedral group in $SO(3)$, and $I^h$ is the 
orientation- reversing extension in $O(3)$.  This  diameter has
also been computed in \cite{M93} and \cite{Gr00}; the details (which are not 
included in \cite{M93} and \cite{Gr00}) are presented here for the sake of 
completeness and because they are important in the calculation of the 
optimal lower bound for quotients of $S^3$ and also give an idea of the 
technique for the diameter of the quotients of $S^2$.   The spherical 
icosahedron has 20 spherical triangular faces, 12 vertices, and 30 edges.
The generators of $I$ are a $\frac{2\pi}{5}$ rotation of a pentagon formed by 
the outer edges of five adjacent faces, a $\frac{2\pi}{3}$ rotation 
about the center 
of a given face, and a $\pi$ rotation about a line through the midpoints of 
opposite edges. The group has order 60. Any triangle can be rotated into 
any other triangle by a combination of rotations, so examine one of these 
triangles.
A rotation of $\frac{2\pi}{3}$ about its center $c$
self-identifies this triangle. The fundamental domain is shaded in  
Figure \ref{Fig:233}, and the diameter is achieved as the length of the spherical 
segment from vertex $v$ to $c$. 
$S^2 / I$ is a triangular shaped inflatable pillow. The space has isotropy 
$C_3$ at $c$,
isotropy $C_2$ at $e$ and isotropy $C_5$ at $v$.
This fundamental domain can be cut in half with a mirror reflection in a 
spherical geodesic 
beginning at $c$ and ending at $e$.  This group $I^h$ is
also a Coxeter group, which can be generated purely by reflections.  The 
fundamental domain
for $S^2/I^h$ is a spherical triangle with vertex angles $\frac{\pi}{2}$ 
at $e$, $\frac{\pi}{3}$ at $c$, and
$\frac{\pi}{5}$ at $v$.    By applying the spherical trigonometry formula 
$\cos A = - \cos B \cos C + \sin A \sin B \cos a$, where $A,B,C$ are vertex 
angles and $a,b,c$ are spherical lengths of opposite edges, one can obtain 
the lower bound.
\end{proof}
\begin{proposition}\label{T:CohomTwoBound}
Let $G$ be a closed subgroup of $O(4)$ acting  on $S^3$ by 
cohomogeneity two.  Then 
\[
\diam(S^3/G)\geq \frac{1}{2}\alpha,
\]
where $\alpha=\arccos (\frac{\tan(\frac{3 \pi}{10})}{\sqrt3})$.
\end{proposition}
\begin{proof}
The only connected group that can act effectively by 
cohomogeneity two on $S^3$ is $T^1=\{e^{i \theta}| \theta \in \rrr \}$.  In the case where $G$ is connected, 
we see that
$T^1$ can act under any of the various possible guises as 
$(z,w) \rightarrow (e^{ik\theta}z, e^{im\theta}w)$.  As such, 
these various group actions are designated by $T_{k,m}$.  Note that this 
action is effective exactly when
$\gcd(k, m)=1$, so we only consider 
these actions.

If $G_0=T^1=T_{1,1}$, then $S^3/G_0$ is isometric to the $2$-sphere of
radius $1/2$, denoted $S^2(1/2)$.  By Lemma \ref{T:componentG}, 
$S^3/G=(S^3/G_0)/(G/G_0)=S^2(1/2)/(G/T_{1,1}).$  
In addition, the action of $G/T_{1,1}$ 
on $S^3/G_0$ is conjugate to the action of a finite subgroup $K$ of $O(3)$ on 
$S^2(1/2)$.  
We obtain the desired lower bound by applying Lemma \ref{T:closedO3} since 
\[
\diam (S^3/G) = \diam (S^2(\frac{1}{2})/K) 
\geq \frac{1}{2}\arccos(\frac{\tan(\frac{3 \pi}{10})}{\sqrt3})
\]
(see also Table \ref{Ta:2Dfin} below).

Now, for the spaces $X_{k,m}=S^3/G_{k,m}$, with $k \neq m$ and $\gcd(k,m)=1$, 
the orbits of $(1,0)$ and $(0,1)$ \emph{cannot}
be interchanged via an isometry. 
This is clear from the previous discussion, since the two exceptional 
orbits are not isometric, having non-isomorphic isotropy subgroups.  
The isotropy subgroup at $(1,0)$ is $C_k$; 
at $(0,1)$, it is $C_m$.   The diameter is realized by a path of length 
$\pi/2$ from $(1,0)$ to $(0,1)$.  Notice that  
$X_{k,m}=S^3/G_{k,m}$ with $k \neq m$ and $\gcd(k,m)=1$ 
is a bad orbifold that is 
topologically a 2-sphere. It has one singular point if  
$k$ or $m$ is 1, and two singular points otherwise.
Any group of isometries of $X_{k,m}$ with $k \neq m$ must fix the orbits of
$(1,0)$ and $(0,1)$, and so quotients of $X_{k,m}$ must always have diameter 
$\pi/2$.

If $G$ is infinite but not connected, we rely on Lemma \ref{T:componentG} 
to reduce the calculation of the diameter of $S^3/G$ to that of 
$S^2(\frac{1}{2})/\Gamma$,
where $\Gamma$ is a finite subgroup of $O(3)$ and then apply 
Lemma \ref{T:closedO3} to obtain the desired lower bound.
\end{proof}
\begin{table}
  \begin{center}
  \caption{Subgroups of $O(3)$ and Corresponding Diameters}\label{Ta:2Dfin}
    \begin{tabular}{lcccccc}
$n$ & \cite{S1891} & \cite{T80} & \cite{Con03} & \cite{Y88} 
& \cite{Int1995} &\text{Diameter} \\
 & $C_n$ & $S^2_{(n,n)}$ & $nn$ & $C_n$ & $n$ & $\pi$\\
odd & $S_{N}$,${C_n}^i$ & $P^2_{(n)}$ & $n\times$ & 
$C_n \cup xC_n$ & $\bar{n}$ & $\pi/2$\\
even & $S_{N}$ & $P^2_{(n)}$ & $n\times$ &  
$C_n\,[\,C_{N}$ & $\bar{N}$& $\pi/2$ \\
odd & ${C_n}^h$ & $D^2_{(n;)}$ & $n*$ &
$C_n\,[\,C_{N}$ & $\bar{N}$& $\pi/2$ \\
even & ${C_n}^h$ & $D^2_{(n;)}$ & $n*$ &  
$C_n \cup xC_n$ & $n$/m& $\pi/2$ \\
odd & ${C_n}^v$ & $D^2_{(;n,n)}$ & $*nn$ &  
$C_n\,[\,D_n$ & $n$m& $\pi$ \\
even & ${C_n}^v$ & $D^2_{(;n,n)}$ & $*nn$ &  
$C_n\,[\,D_n$ & $n$mm & $\pi$ \\
   &&&&&\\
odd & $D_n$ & $S^2_{(2,2,n)}$ & $22n$ &  $D_n$ & 
$n$2  & $\pi/2$\\
even & $D_n$ & $S^2_{(2,2,n)}$ & $22n$ &  $D_n$ & 
$n$22& $\pi/2$ \\
odd & ${D_n}^h$ & $D^2_{(;2,2,n)}$ & $*22n$ & 
$D_n\,[\,D_{N}$ & $\bar{N}$2m& $\pi/2$ \\
even & ${D_n}^h$ & $D^2_{(;2,2,n)}$ & $*22n$ &
$D_n \cup xD_n$ & $n$/mmm & $\pi/2$\\
odd & ${D_n}^d$ & $D^2_{(2;n)}$ & $2{*}n$ &  
$D_n \cup xD_n$ & $\bar{n}$m& $\pi/2$ \\
even & ${D_n}^d$ & $D^2_{(2;n)}$ & $2{*}n$ & 
$D_n\,[\,D_{N}$ & $\bar{N}$2m & $\pi/2$\\
   &&&&&\\ & $T$ & $S^2_{(2,3,3)}$ & $332$ &  $T$ & $23$ &
$\arccos{\frac{1}{3}}$\\
 & $T^d$ & $D^2_{(;2,3,3)}$ & $*332$ &$T\,[\,O$ & 
$\bar{4}$3m &$\arccos{\frac{1}{3}}$\\
 & $T^h$ & $D^2_{(3;2)}$ & $3{*}2$ &  $T \cup xT$ & 
m$\bar{3}$ &$\arccos{\frac{1}{\sqrt{3}}}$\\
   &&&&&\\
 & $O$ & $S^2_{(2,3,4)}$ & $432$ &  $O$ & 
$432$&$\arccos{\frac{1}{\sqrt{3}}}$ \\
 & $O^h$ & $D^2_{(;2,3,4)}$ & $*432$ &  
$O \cup xO$ & m$\bar{3}$m &$\arccos{\frac{1}{\sqrt{3}}}$\\
   &&&&&\\
 & $I$ & $S^2_{(2,3,5)}$ & $532$ &  $I$ & $235$ &
$\arccos{\frac{\tan(\frac{3 \pi}{10})}{\sqrt3}}$\\
 & $I^h$ & $D^2_{(;2,3,5)}$ & $*532$ &  $I \cup xI$ & 
m$\bar{3}\bar{5}$ &$\arccos{\frac{\tan(\frac{3 \pi}{10})}{\sqrt3}}$\\
    \end{tabular}
  \end{center}
\end{table}
In Table \ref{Ta:2Dfin}, the second \cite{S1891} column describes 
the geometry of the group actions.  
The third \cite{T80} and fourth \cite{Con03} columns list 
the resulting topological orbit spaces (see pictures in \cite{M85}).
The fifth column \cite{Y88}, is given for better understanding of the algebraic
structure of the groups; both this
and the last column \cite{Int1995} focus on inversions instead of  
reflections.
For the ``$C$'' groups, $n \ge 1$; for the ``$D$'' groups, $n \ge 2$.  The 
special case ``$2$/mmm'' is
written ``mmm''.  ``$N$'' stands for the number $2n$. ``$x$'' stands for the 
antipodal map.
When $H$ is a subgroup of index two of $G \subseteq SO(3)$, 
define $H\ [\ G := H \cup x(G-H)$.

\section{Cohomogeneity Three}\label{S:cohomThree}

In this section we compute a lower bound on the diameters of $S^3/G$ where 
$G$ is finite in $O(4)$.

\subsection{Classification of finite subgroups of $\mathbf{O(4)}$}\label{S: classifO4}

While this discussion will basically follow the 
treatment in Du\,Val \cite{D64}, the reader
should note that Threlfall and Seifert \cite{TS30,TS32} 
classify finite subgroups of $SO(4)$.  
Conway and Smith \cite[Chapter 4]{Con03} also have a 
classification of subgroups of $O(4)$.

The central idea in the classification of finite subgroups of $O(4)$, 
up to conjugacy,
is that $SO(4)$ is ``almost a product''.  More precisely, there are 
2-to-1 homomorphisms
$S^3 \times S^3 \longrightarrow SO(4) \longrightarrow SO(3) \times SO(3)$ 
(the former homomorphism was called $\eta$ in section \ref{S:intro}). 

They are defined by thinking of $S^3$ as the set of unit quaternions:
\[
(\mathbf{p}_1,\mathbf{p}_2) \overset{\eta}{\mapsto}
(\mathbf{q} \mapsto \mathbf{p}_1 \mathbf{q} \mathbf{p}_2^{-1}) \mapsto
((\tilde{\mathbf{q}} \mapsto \mathbf{p}_1 \tilde{\mathbf{q}}  
\mathbf{p}_1^{-1}),
  (\tilde{\mathbf{q}} \mapsto \mathbf{p}_2 \tilde{\mathbf{q}}  
\mathbf{p}_2^{-1})),
\]
where $\mathbf{q} := q_1+ q_2\mathbf{i} + q_3\mathbf{j} + q_4\mathbf{k}$
and $\tilde{\mathbf{q}} := q_2\mathbf{i} + q_3\mathbf{j} +  
q_4\mathbf{k}$.

Finite subgroups of 
$SO(4)$ are classified by combining the well-known classification of finite 
subgroups
of $SO(3)$ (cyclic, dihedral, tetrahedral, octahedral, icosahedral ---
see, e.g., \cite{Y88})
and the less-well-known, but elementary classification of subgroups of 
product groups 
(sketched below, but see also \cite[pages 63--64]{Ha59}).

If $G$ denotes a finite subgroup of $SO(4)$, then let $\hat G$ denote  
its
inverse image $\eta^{-1}(G)$ in $S^3 \times S^3$.  We define the following subgroups  
of $S^3$:
\begin{align}
\mathbf{L} &:= \{\ell : (\ell,r) \in \hat G \text{ for some } r \} \notag\\
\mathbf{R} &:= \{r : (\ell,r) \in \hat G \text{ for some } \ell \} \notag\\
\textsc{l} &:= \{\ell : (\ell,1) \in \hat G \} \notag\\
\textsc{r} &:= \{r: (1,r) \in \hat G \} \notag
\end{align}
It can be shown that 
\[
\textsc{l} =
\ker (\lambda':\mathbf{L} \rightarrow \hat G/(\textsc{l}\times\textsc{r}) )
\mbox{ and }
 \textsc{r} =
\ker (\rho':\mathbf{R} \rightarrow \hat G/(\textsc{l}\times\textsc{r}) ),
\]
inducing isomorphisms
$\lambda: \mathbf{L}/\textsc{l} \rightarrow  \hat G/(\textsc{l}\times\textsc{r})$
and
$\rho:\mathbf{R}/\textsc{r} \rightarrow  \hat G/(\textsc{l}\times\textsc{r})$
which, when composed back-to-back, give
an isomorphism $\phi = \rho^{-1} \circ \lambda$ from
$\mathbf{L}/\textsc{l}$ to $\mathbf{R}/\textsc{r}$.
The group $G$ is denoted
$(\mathbf{L}/\textsc{l};\mathbf{R}/\textsc{r};\phi)$.  
Often $\phi$ is omitted if the isomorphism is 
``obvious''; compare \cite[page 54]{D64}.   
The order of $G$ will equal 
$[|\mathbf{R}||\mathbf{L}|/(|\mathbf{L}/\textsc{l}|)]/2 
= |\mathbf{R}||\textsc{l}|/2$.

The only possibilities for $\mathbf{L}$ and $\mathbf{R}$ are the finite 
subgroups of $S^3$ which
are inverse images under the 2-to-1 homomorphism $S^3 \longrightarrow SO(3)$ 
of finite subgroups of 
$SO(3)$, or in other words, finite subgroups of $S^3$ containing the kernel 
$\{\pm 1\}$ of
that homomorphism.  Hence they are conjugate to (exactly) one of the 
following ``binary''
groups:
\begin{align}
\mathbf{C}_{2n} &:= 
\{ \cos(2m\pi/2n) + \sin(2m\pi/2n)\mathbf{k}: m =  0,1,\dots,2n-1 \} \notag\\
& \quad (n \ge 1) \notag\\
\mathbf{D}_n &:= \mathbf{C}_{2n} \cup\{ \cos(2m\pi/2n)\mathbf{i} + 
\sin(2m\pi/2n)\mathbf{j}: m =  0,1,\dots,2n-1 \} \notag\\
& \quad (n \ge 2) \notag\\
\mathbf{T} &:= \mathbf{D}_2 \cup 
\{ \pm \frac{1}{2}  \pm \frac{1}{2}\mathbf{i}  \pm \frac{1}{2}\mathbf{j} 
\pm \frac{1}{2}\mathbf{k} \}\notag\\
\mathbf{O} &:= \mathbf{T} \cup
\{\pm 1/\sqrt{2} \pm (1/\sqrt{2})\mathbf{i} \} \cup
\{\pm 1/\sqrt{2} \pm (1/\sqrt{2})\mathbf{j} \}\notag\\
&\cup \{\pm 1/\sqrt{2} \pm (1/\sqrt{2})\mathbf{k} \} \cup
 \{\pm (1/\sqrt{2})\mathbf{i} \pm (1/\sqrt{2})\mathbf{j} \}\notag\\
 &\cup
\{\pm (1/\sqrt{2})\mathbf{i} \pm (1/\sqrt{2})\mathbf{k} \} \cup
\{\pm (1/\sqrt{2})\mathbf{j} \pm (1/\sqrt{2})\mathbf{k} \}\notag\\
\mathbf{I} &:= \mathbf{T} \cup (1/2)((\tau - 1) + \tau \mathbf{i} + 
\mathbf{j})\mathbf{T}
 \cup (1/2)(-\tau + \mathbf{i} + (\tau - 1)\mathbf{j})\mathbf{T}\notag\\
& \cup (1/2)(-\tau - \mathbf{i} + (1 - \tau)\mathbf{j})\mathbf{T}
 \cup (1/2)((\tau - 1) - \tau \mathbf{i} - \mathbf{j})\mathbf{T}\notag,
\end{align}
where $\tau := (\sqrt{5}+1)/2$.
These groups have orders $2n$, $4n$, $24$, $48$, and $120$, respectively.  
It may be worth noting that the given coset representatives for $\mathbf{I}$ 
form
a cyclic group of order five.  If the signs of all the coefficients of 
$\sqrt{5}$
in elements of $\mathbf{I}$ are reversed, then a group $\mathbf{I}^\dag$ is
obtained, which is conjugate to $\mathbf{I}$ in $S^3$, and  has the
property that 
$\mathbf{I} \cap \mathbf{I}^\dag = \mathbf{T}$ \cite[page 55]{D64}.   
``Sign-reversal'' defines an isomorphism 
$\phi_\dag: \mathbf{I} \rightarrow \mathbf{I}^\dag$,
whose inverse is also accomplished by sign-reversal (i.e., by
a different sort of ``conjugation'', in the field 
$\mathbb{Q}(\sqrt{5})$).

We arrive at 41 families of finite subgroups of $SO(4)$, 33 of which contain
the central element (-1 times the identity matrix) and which therefore
equal the inverse image of their projections to $SO(3) \times SO(3)$.  
The numbering convention follows Du\,Val \cite{D64} and goes back
to Goursat \cite{G1889}, 
who classified the finite subgroups of 
$\Isom(\mathbb{RP}^3) \cong SO(3) \times SO(3)$, 
though in some places we are forced to interpolate extra families 
to cover gaps in that enumeration. 
The first 33 families are listed in Table \ref{Ta:finI}, where $m,n,r \ge 1$, 
$\gcd(s,r)=1$, and $0 \le s < r/2$
\cite[page 55]{D64}.
Furthermore, $\phi_s: \mathbf{C}_{2mr}/\mathbf{C}_{2m}\rightarrow  
\mathbf{C}_{2nr}/\mathbf{C}_{2n}$
is the isomorphism which takes the coset
$(\cos(2\pi/2mr) + \sin(2\pi/2mr)\mathbf{k}) \mathbf{C}_{2m}$
to the coset $(\cos(2s\pi/2nr) + \sin(2s\pi/2nr)\mathbf{k})  
\mathbf{C}_{2n}$.
Similarly, $\psi_s: \mathbf{D}_{mr}/\mathbf{C}_{2m} \rightarrow  
\mathbf{D}_{nr}/\mathbf{C}_{2n}$
is the isomorphism mapping cosets with representatives in  
$\mathbf{C}_{2mr}$ as above,
while taking the coset $\mathbf{i} \mathbf{C}_{2m}$ to the
coset $\mathbf{i} \mathbf{C}_{2n}$.
Finally, $\phi_\dag : \mathbf{I} \rightarrow \mathbf{I}^\dag$ induces  
an isomorphism $\tilde{\phi}_\dag : \mathbf{I}/\mathbf{C}_2 \rightarrow  
\mathbf{I}^\dag/\mathbf{C}_2$.
\begin{table}
  \begin{center}
  \caption{Finite Subgroups of $O(4)$, part I}\label{Ta:finI}
    \[\begin{array}{lclc} 
1. & (\mathbf{C}_{2mr}/\mathbf{C}_{2m};\mathbf{C}_{2nr}/\mathbf{C}_{2n};
\phi_s)  &
17. & (\mathbf{D}_{2m}/\mathbf{D}_m;\mathbf{O}/\mathbf{T}) \\
2. & (\mathbf{C}_{2m}/\mathbf{C}_{2m}; \mathbf{D}_n/\mathbf{D}_n)  &
18. & (\mathbf{D}_{6m}/\mathbf{C}_{2m}; \mathbf{O}/\mathbf{D}_2) \\
3. & (\mathbf{C}_{4m}/\mathbf{C}_{2m}; \mathbf{D}_{n}/\mathbf{C}_{2n})  &
19. & (\mathbf{D}_m/\mathbf{D}_m; \mathbf{I}/ \mathbf{I}) \\
4. & (\mathbf{C}_{4m}/\mathbf{C}_{2m}; \mathbf{D}_{2n}/\mathbf{D}_n) & 
20. & (\mathbf{T}/\mathbf{T}; \mathbf{T}/ \mathbf{T}) \\
5. & (\mathbf{C}_{2m}/\mathbf{C}_{2m}; \mathbf{T}/\mathbf{T}) & 
21. & (\mathbf{T}/\mathbf{C}_2; \mathbf{T}/ \mathbf{C}_2) \\
6. & (\mathbf{C}_{6m}/\mathbf{C}_{2m}; \mathbf{T}/\mathbf{D}_2) & 
22. & (\mathbf{T}/\mathbf{D}_2; \mathbf{T}/ \mathbf{D}_2) \\
7. & (\mathbf{C}_{2m}/\mathbf{C}_{2m}; \mathbf{O}/\mathbf{O}) & 
23. & (\mathbf{T}/\mathbf{T}; \mathbf{O}/ \mathbf{O}) \\
8. & (\mathbf{C}_{4m}/\mathbf{C}_{2m}; \mathbf{O}/\mathbf{T}) & 
24. & (\mathbf{T}/\mathbf{T}; \mathbf{I}/ \mathbf{I}) \\
9. & (\mathbf{C}_{2m}/\mathbf{C}_{2m}; \mathbf{I}/\mathbf{I}) & 
25. & (\mathbf{O}/\mathbf{O}; \mathbf{O}/ \mathbf{O}) \\
10. & (\mathbf{D}_m/\mathbf{D}_m; \mathbf{D}_n/\mathbf{D}_n) & 
26. & (\mathbf{O}/\mathbf{C}_2; \mathbf{O}/ \mathbf{C}_2) \\
11. &  
(\mathbf{D}_{mr}/\mathbf{C}_{2m};\mathbf{D}_{nr}/\mathbf{C}_{2n}; 
\psi_s) &
27. & (\mathbf{O}/\mathbf{D}_2; \mathbf{O}/ \mathbf{D}_2) \\
11a. &
(\mathbf{D}_{2m}/\mathbf{C}_{2m};\mathbf{D}_{2n}/\mathbf{C}_{2n}; 
\psi_\#) & 
28. & (\mathbf{O}/\mathbf{T}; \mathbf{O}/ \mathbf{T}) \\
12. & (\mathbf{D}_{2m}/\mathbf{D}_m; \mathbf{D}_{2n}/\mathbf{D}_n) & 
29. & (\mathbf{O}/\mathbf{O}; \mathbf{I}/ \mathbf{I}) \\
13. & (\mathbf{D}_{2m}/\mathbf{D}_m; \mathbf{D}_{n}/\mathbf{C}_{2n}) & 
30. & (\mathbf{I}/\mathbf{I}; \mathbf{I}/ \mathbf{I}) \\
14. & (\mathbf{D}_m/\mathbf{D}_m; \mathbf{T}/\mathbf{T}) & 
31. & (\mathbf{I}/\mathbf{C}_2;\mathbf{I}/\mathbf{C}_2) \\
15. & (\mathbf{D}_m/\mathbf{D}_m; \mathbf{O}/\mathbf{O}) &
32. & (\mathbf{I}^\dag/\mathbf{C}_2;\mathbf{I}/\mathbf{C}_2;
\tilde{\phi}_\dag^{-1}) \\
16. & (\mathbf{D}_m/\mathbf{C}_{2m};\mathbf{O}/\mathbf{T}) & & \\
    \end{array}\]
  \end{center}
\end{table}
As noted in \cite[page 585]{TS32} and \cite[page 50]{Con03}, Goursat
and Du\,Val
omit a family of the form ($m,n \ge 2$)
\[
11a. \quad  
(\mathbf{D}_{2m}/\mathbf{C}_{2m};\mathbf{D}_{2n}/\mathbf{C}_{2n}; 
\psi_\#)
\]
where the common quotient group is the Klein four-group, and $\psi_\#$  
is the isomorphism which takes the coset
$(\cos(\pi/m) + \sin(\pi/m)\mathbf{k}) \mathbf{C}_{2m}$
to the coset
$\mathbf{i} \mathbf{C}_{2n}$
and conversely takes the coset
$\mathbf{i} \mathbf{C}_{2m}$
to the coset
$(\cos(\pi/n) + \sin(\pi/n)\mathbf{k}) \mathbf{C}_{2n}$
(which does not respect the cyclic subgroups of index two in these
binary dihedral groups).

There are also 8 families of subgroups of $SO(4)$ which do not contain the 
central element.  
It follows that \textsc{l} and \textsc{r} for such groups must be cyclic of 
odd order, 
since each of the other subgroups of $S^3$ contains the quaternion $-1$;
we therefore extend the notation $\mathbf{C}_{2n}$ to allow odd 
subscripts.     
In addition to the conditions on $m,n,r,s$ given above, in Table \ref{Ta:finII},
$m$ and $n$ are both odd integers.
\begin{table}
  \begin{center}
  \caption{Finite Subgroups of $O(4)$, part II}\label{Ta:finII}
    \[\begin{array}{lclc}
1'. \quad & 
(\mathbf{C}_{2mr}/\mathbf{C}_m;\mathbf{C}_{2nr}/\mathbf{C}_n;\phi_s) &
26'. \quad & (\mathbf{O}/\mathbf{C}_1;\mathbf{O}/\mathbf{C}_1;id)\\
11'. \quad & 
(\mathbf{D}_{mr}/\mathbf{C}_m;\mathbf{D}_{nr}/\mathbf{C}_n;\psi_s) &
26''. \quad &(\mathbf{O}/\mathbf{C}_1;\mathbf{O}/\mathbf{C}_1;\xi)\\
11a'. \quad & 
(\mathbf{D}_{2m}/\mathbf{C}_{m};\mathbf{D}_{2n}/\mathbf{C}_{n};\psi_\#) 
&
31'. \quad & (\mathbf{I}/\mathbf{C}_1;\mathbf{I}/\mathbf{C}_1) \\
21'. \quad & (\mathbf{T}/\mathbf{C}_1;\mathbf{T}/\mathbf{C}_1)&
32'. \quad & 
(\mathbf{I}^\dag/\mathbf{C}_1;\mathbf{I}/\mathbf{C}_1;\phi_\dag^{-1})
    \end{array}\]
  \end{center}
\end{table}
The automorphism in \#$26''$, 
$\xi: \mathbf{O} \rightarrow \mathbf{O}$, is the identity on 
$\mathbf{T}$,
and multiplies all other elements by $-1$.  It cannot be induced by 
conjugation in $S^3$, and 
hence groups \#$26'$ and \#$26''$ are not conjugate in $SO(4)$. 

It remains to list, up to conjugacy, the finite subgroups of $O(4)$ which 
contain orientation-reversing transformations.  
We can write an arbitrary
element of $O(4) - SO(4)$ as the composition of the
particular orientation-reversing map
$\mathbf{q} \mapsto \overline{\mathbf{q}}$
(linear, mapping $1\mapsto 1$, $\mathbf{i}\mapsto -\mathbf{i}$,
$\mathbf{j}\mapsto -\mathbf{j}$, and $\mathbf{k}\mapsto -\mathbf{k}$),
followed by an arbitrary orientation-preserving map
$\mathbf{q} \mapsto \mathbf{a}\mathbf{q}\mathbf{b}$, hence in the form 
$\mathbf{q} \mapsto \mathbf{a} \overline{\mathbf{q}} \mathbf{b}$
(where $\mathbf{a},\mathbf{b} \in S^3$).  This representation is unique
up to multiplying both $\mathbf{a}$ and $\mathbf{b}$ by $-1$.  

It follows from the identity 
$\overline{\mathbf{q}_1\mathbf{q}_2} = \overline{\mathbf{q}_2}\ 
\overline{\mathbf{q}_1}$ 
that conjugation by the orientation-reversing map 
$\mathbf{q} \mapsto \overline{\mathbf{q}}$
takes an element of $SO(4)$ covered by $(\ell,r) \in S^3 \times S^3$ 
to one covered by $(r,\ell)$.  
Hence, more generally, for a finite subgroup $G$ of $O(4)$ containing
orientation-reversing elements, the groups $\mathbf{L}$ and $\mathbf{R}$
describing $G \cap SO(4)$ must be conjugate.  Indeed, $\mathbf{L} = 
\mathbf{R}$,
except when $G \cap SO(4)$ equals group \#$32$ or group \#$32'$.

Following Du\,Val's classification, we start with those subgroups of $G$ of 
$O(4)$ which contain the central element (and a few which do not, namely families  
\#$33$, $35$, and $36$, when $n$ is odd).  
Du\,Val's notation adds a superscript asterisk to the symbol for the 
orientation-preserving
subgroup, in order to indicate the presence 
of orientation-reversing elements.  In most cases, we specify the extension 
by describing
$\{(\mathbf{a},\mathbf{b}):
\mathbf{q} \mapsto \mathbf{a} \overline{\mathbf{q}} \mathbf{b} \in G
\}$.
The basic conditions on the integers $n,r,s,h,k$ in Table \ref{Ta:finIII} are
$n,r \ge 1$, $0 \le s,h,k <r$, $\gcd(s,r)=1$ and $rn$ is even;
further conditions are as follows:

\smallskip

[1]:
$s^2 \equiv 1$, $h(s-1)\equiv 0 \pmod{r}$; extending element    
$\mathbf{q} \mapsto \mathbf{p}^h \overline{\mathbf{q}}$, where
$\mathbf{p} := \cos\frac{2\pi}{nr}+\sin\frac{2\pi}{nr}\mathbf{k}$.

[2]:
$s^2 \equiv 1\pmod{r}$, $h\equiv k \pmod{2}$,  
$(h-k)(s-1)\equiv (h+k)(s+1) \equiv 0 \pmod{2r}$; extending element 
$\mathbf{q} \mapsto
\mathbf{p}^{\frac{1}{2}h}\overline{\mathbf{q}}\mathbf{p}^{\frac{1}{2}k}$,
using $\mathbf{p}^{\frac{1}{2}}$ to denote
$\cos\frac{\pi}{nr}+\sin\frac{\pi}{nr}\mathbf{k}$.

[3]:
$s^2+1 \equiv 0\pmod{r}$, $h\equiv k \pmod{2}$, 
$h+k \equiv s(h-k), k-h \equiv s(k+k) \pmod{2r}$; extending element
$\mathbf{q} \mapsto
\mathbf{i}\mathbf{p}^{\frac{1}{2}h}\overline{\mathbf{q}}
\mathbf{p}^{\frac{1}{2}k}$, with $\mathbf{p}^{\frac{1}{2}}$ as in note [2].

[4]:
$\mathbf{a} = \mathbf{p}^\dag \mathbf{t}'$, $\mathbf{b} =
\pm(\mathbf{p} \mathbf{t}')^{-1}$, with $\mathbf{p} \in \mathbf{I}$, and 
$\mathbf
{t}' \in \mathbf{O} - \mathbf{T}$. It suffices to let $\mathbf{t}'$ be any fixed element of 
$\mathbf{O} - \mathbf{T}$, such as $(1/\sqrt{2})\mathbf{i} + (1/\sqrt{2})\mathbf{j}$.

\smallskip

\begin{table}
  \begin{center}
  \caption{Finite Subgroups of $O(4)$, part III}\label{Ta:finIII}
  \[\begin{array}{lcc} 
33. &
(\mathbf{C}_{nr}/\mathbf{C}_{n};\mathbf{C}_{nr}/\mathbf{C}_{n};
\phi_s)_h^*;& [1] \\
34. & (\mathbf{D}_n/\mathbf{D}_n; \mathbf{D}_n/\mathbf{D}_n)^*;&   
\mathbf{a},\mathbf{b} \in \mathbf{D}_n  \\
35. & (\mathbf{D}_{\frac{1}{2}nr}/\mathbf{C}_{n};
\mathbf{D}_{\frac{1}{2}nr}/\mathbf{C}_{n};\phi_s)_{h,k}^*;& [2]  \\
35a. & (\mathbf{D}_{2n}/\mathbf{C}_{2n};
\mathbf{D}_{2n}/\mathbf{C}_{2n};\psi_\#)^*;&
\mathbf{a},\mathbf{b}\in \mathbf{D}_{2n},
\mathbf{b}\mathbf{C}_{2n} = \psi_\#(\mathbf{a}\mathbf{C}_{2n})  \\
36. & (\mathbf{D}_{\frac{1}{2}nr}/\mathbf{C}_{n};
\mathbf{D}_{\frac{1}{2}nr}/\mathbf{C}_{n};\phi_s)_{h,k-}^*;& [3] \\
37. & (\mathbf{D}_{2n}/\mathbf{D}_n; \mathbf{D}_{2n}/\mathbf{D}_n)^*;& 
\mathbf{a},\mathbf{b} \in \mathbf{D}_{2n}, \mathbf{a}\mathbf{D}_n = 
\mathbf{b}\mathbf{D}_n \\ 
38. & (\mathbf{D}_{2n}/\mathbf{D}_n; \mathbf{D}_{2n}/\mathbf{D}_n)^*_-;& 
\mathbf{a},\mathbf{b} \in \mathbf{D}_{2n}, \mathbf{a}\mathbf{D}_n \ne 
\mathbf{b}\mathbf{D}_n \\ 
39. & (\mathbf{T}/\mathbf{C}_2; \mathbf{T}/ \mathbf{C}_2)^*_c;& \mathbf{a} 
\in \mathbf{T}, \mathbf{b} = \pm\mathbf{a}^{-1} \\ 
40. & (\mathbf{T}/\mathbf{C}_2; \mathbf{T}/ \mathbf{C}_2)^*;& \mathbf{a} \in 
\mathbf{O}-\mathbf{T}, \mathbf{b} = \pm\mathbf{a}^{-1} \\ 
41. & (\mathbf{T}/\mathbf{D}_2; \mathbf{T}/ \mathbf{D}_2)^*;& \mathbf{a},
\mathbf{b} \in \mathbf{T}, \mathbf{ab} \in \mathbf{D}_2 \\
42. & (\mathbf{T}/\mathbf{D}_2; \mathbf{T}/ \mathbf{D}_2)^*_-;& \mathbf{a},
\mathbf{b} \in \mathbf{O}-\mathbf{T}, \mathbf{ab}^{-1} \in \mathbf{D}_2 \\
43. & (\mathbf{T}/\mathbf{T}; \mathbf{T}/ \mathbf{T})^*;& \mathbf{a},
\mathbf{b} \in \mathbf{T} \\
44. & (\mathbf{O}/\mathbf{C}_2; \mathbf{O}/ \mathbf{C}_2)^*;& \mathbf{a} 
\in \mathbf{O}, \mathbf{b} = \pm\mathbf{a}^{-1} \\
45. & (\mathbf{O}/\mathbf{T}; \mathbf{O}/ \mathbf{T})^*;& \mathbf{a},
\mathbf{b} \in \mathbf{O}, \mathbf{aT} = \mathbf{bT} \\
46. & (\mathbf{O}/\mathbf{T}; \mathbf{O}/ \mathbf{T})^*_-;& \mathbf{a},
\mathbf{b} \in \mathbf{O}, \mathbf{aT} \ne \mathbf{bT} \\
47. & (\mathbf{O}/\mathbf{D}_2; \mathbf{O}/ \mathbf{D}_2)^*;& \mathbf{a},
\mathbf{b} \in \mathbf{O}, \mathbf{ab} \in \mathbf{D}_2 \\
48. & (\mathbf{O}/\mathbf{O}; \mathbf{O}/ \mathbf{O})^*;& \mathbf{a},
\mathbf{b} \in \mathbf{O} \\
49. & (\mathbf{I}/\mathbf{C}_2; \mathbf{I}/ \mathbf{C}_2)^*;& \mathbf{a} 
\in \mathbf{I}, \mathbf{b} = \pm\mathbf{a}^{-1} \\
50. & (\mathbf{I}/\mathbf{I}; \mathbf{I}/ \mathbf{I})^*;& \mathbf{a},
\mathbf{b} \in \mathbf{I} \\
51. & (\mathbf{I}^\dag/\mathbf{C}_2;\mathbf{I}/\mathbf{C}_2;
\tilde{\phi}_\dag^{-1})^*;& [4]
    \end{array}\]
  \end{center}
\end{table}
In addition, a few of these groups have subgroups of index two which contain
orientation-reversing elements, but do not contain the central element.
These are listed in Table \ref{Ta:finIV}.  
Groups \#$44pm$ and \#$44mp$ do not appear in \cite[page 61]{D64}, but
are listed as
$\pm\frac{1}{24}[O \times \overline{O}]\cdot 2_3$
and $\pm\frac{1}{24}[O \times \overline{O}]\cdot 2_1$, respectively,
in \cite[page 47]{Con03}.  Extra conditions on the groups in this table are
as follows:

\smallskip

[5]:
($\mathbf{a},\mathbf{b}\in \mathbf{C}_n$) or
($\mathbf{a}\in -\mathbf{k}\mathbf{C}_{n},
\mathbf{b}\in \mathbf{i}\mathbf{C}_{n}$) or
($\mathbf{a}\in \mathbf{i}\mathbf{C}_{n},
\mathbf{b}\in -\mathbf{k}\mathbf{C}_{n}$) or
($\mathbf{a},\mathbf{b}\in \mathbf{j}\mathbf{C}_n$).

[6]:
($\mathbf{a},-\mathbf{b}\in \mathbf{C}_n$) or
($\mathbf{a}\in \mathbf{k}\mathbf{C}_{n},
\mathbf{b}\in \mathbf{i}\mathbf{C}_{n}$) or
($\mathbf{a}\in \mathbf{i}\mathbf{C}_{n},
\mathbf{b}\in \mathbf{k}\mathbf{C}_{n}$) or
($\mathbf{a},-\mathbf{b}\in \mathbf{j}\mathbf{C}_n$).

[7]: $\mathbf{a} = \mathbf{p}^\dag \mathbf{t}'$, $\mathbf{b} = 
(\mathbf{p} \mathbf{t}')^{-1}$,
with $\mathbf{p} \in \mathbf{I}$, and $\mathbf{t}' \in \mathbf{O} - 
\mathbf{T}$.

[8]: $\mathbf{a} = \mathbf{p}^\dag \mathbf{t}'$, $\mathbf{b} = 
-(\mathbf{p} \mathbf{t}')^{-1}$,
with $\mathbf{p} \in \mathbf{I}$, and $\mathbf{t}' \in \mathbf{O} - 
\mathbf{T}$.

\smallskip

\begin{table}
  \begin{center}
  \caption{Finite Subgroups of $O(4)$, part IV}\label{Ta:finIV}
     \[\begin{array}{lcc}
35ap. & (\mathbf{D}_{2n}/\mathbf{C}_{n};
\mathbf{D}_{2n}/\mathbf{C}_{n};\psi_\#)^*;& [5]  \\
35am. & (\mathbf{D}_{2n}/\mathbf{C}_{n};
\mathbf{D}_{2n}/\mathbf{C}_{n};\psi_\#)^*_{-};& [6]  \\
39p. & (\mathbf{T}/\mathbf{C}_1; \mathbf{T}/ \mathbf{C}_1)^*_c;& \mathbf{a} 
\in \mathbf{T}, \mathbf{b} = \mathbf{a}^{-1} \\
39m. & (\mathbf{T}/\mathbf{C}_1; \mathbf{T}/ \mathbf{C}_1)^*_{c-};& 
\mathbf{a} \in \mathbf{T}, \mathbf{b} = -\mathbf{a}^{-1} \\
40p. & (\mathbf{T}/\mathbf{C}_1; \mathbf{T}/ \mathbf{C}_1)^*;& \mathbf{a} 
\in \mathbf{O}-\mathbf{T}, \mathbf{b} = \mathbf{a}^{-1} \\ 
40m. & (\mathbf{T}/\mathbf{C}_1; \mathbf{T}/ \mathbf{C}_1)^*_{-};& \mathbf{a} 
\in \mathbf{O}-\mathbf{T}, \mathbf{b} = -\mathbf{a}^{-1} \\ 
44p. & (\mathbf{O}/\mathbf{C}_1; \mathbf{O}/ \mathbf{C}_1;id)^*;& \mathbf{a} 
\in \mathbf{O}, \mathbf{b} = \mathbf{a}^{-1} \\
44m. & (\mathbf{O}/\mathbf{C}_1; \mathbf{O}/ \mathbf{C}_1;id)^*_{-};& 
\mathbf{a} \in \mathbf{O}, \mathbf{b} = -\mathbf{a}^{-1} \\
44pm. & (\mathbf{O}/\mathbf{C}_1; \mathbf{O}/ \mathbf{C}_1;\xi)^*_{+-};& 
\mathbf{b}^{-1} = \mathbf{a} \in \mathbf{T} \mbox{ or }
{-\mathbf{b}^{-1}} = \mathbf{a} \in \mathbf{O}-\mathbf{T} \\
44mp. & (\mathbf{O}/\mathbf{C}_1; \mathbf{O}/ \mathbf{C}_1;\xi)^*_{-+};& 
-\mathbf{b}^{-1} = \mathbf{a} \in \mathbf{T} \mbox{ or }
\mathbf{b}^{-1} = \mathbf{a} \in \mathbf{O}-\mathbf{T} \\
49p. & (\mathbf{I}/\mathbf{C}_1; \mathbf{I}/ \mathbf{C}_1)^*;& \mathbf{a} 
\in \mathbf{I}, \mathbf{b} = \mathbf{a}^{-1} \\
49m. & (\mathbf{I}/\mathbf{C}_1; \mathbf{I}/ \mathbf{C}_1)^*_{-};& \mathbf{a} 
\in \mathbf{I}, \mathbf{b} = -\mathbf{a}^{-1} \\
51p. & (\mathbf{I}^\dag/\mathbf{C}_1;\mathbf{I}/\mathbf{C}_1;
\tilde{\phi}_\dag^{-1})^*;& [7] \\
51m. & (\mathbf{I}^\dag/\mathbf{C}_1;\mathbf{I}/\mathbf{C}_1;
\tilde{\phi}_\dag^{-1})^*_{-};& [8]
    \end{array}\]
  \end{center}
\end{table}
For each of these subgroups $G$ of $O(4)$, we will obtain a lower
bound on the diameter of the orbifold $S^3/G$ by
maximizing the distance from the orbit of
the quaternion $1$, under the group action, to another orbit.   
To that end, we'll define the {\it pre-fundamental domain} 
of a finite group action on $S^3$ to be 
the intersection of half-spheres formed by the set
of points which are closer to the quaternion $1$ than
to any other image of $1$ (equivalent to the Voronoi cell of $1$
with respect to its $G$-orbit; when, in addition, $1$ is fixed only by 
the identity element of $G$, it is a Dirichlet domain). 
The geodesic segment from $1$ to
any point in the pre-fundamental domain realizes the distance
in the orbifold between the equivalence classes represented
by those points.  So the distance from $1$ to the farthest vertex 
of the pre-fundamental domain gives a lower bound for the
diameter of the orbifold.  We expect that this lower bound
will be sharp in most cases, based on the fact that $1$ has a
large isotropy subgroup under the action of $G$ (roughly speaking,
because $\mathbf{L}$ and $\mathbf{R}$ will have large intersection, compared
to other conjugates in $S^3$).  Hence, in the orbifold $S^3/G$, the
point corresponding to $1$ is at the vertex of a ``sharp'' cone, which
should tend to ``push it away'' from ``the rest of the orbifold''.
This assertion is motivated by the fact that singular points 
in a hyperbolic $3$-orbifold which are locally modeled on $\mathbb{H}^3/C$ 
(where $C$ is a large cyclic group acting by rotation) 
are contained in the middle of fat Margulis tubes, hence are far from
the ``thick part'' of the orbifold \cite{Me88}.  This intuition is quite rough,
but explains our choice of $1$ for one end of the geodesic segment 
which realizes our lower bound for the diameter.  In some cases, noted below, we
can prove that our bound is sharp.  
The coordinates of the vertices of a pre-fundamental domain can
be calculated by linear algebra once triples of points in the orbit of $1$ are
found which are both close to $1$ and close to each other;
there are three linear constraints since the vertex must be
equidistant from $1$ and each point in the triple, and also the
vertex must have unit length (and make an acute angle with $1$).
We used Maple${}^\mathrm{TM}$ software to handle the messier situations \cite{M05}.
See \cite{Du94} for more on pre-fundamental domains;
in particular, a fundamental domain for $S^3/G$ can be obtained
by intersecting the pre-fundamental domain of $G$ with a cone which is
a fundamental domain for the subgroup of $G$ which fixes $1$.  

\subsection{Diameters for the fibering subgroups of $\mathbf{SO(4)}$}\label{S:fibSO}

The subgroups $G$ of $SO(4)$ for which the corresponding orbifolds
$S^3/G$ admit a fibering over a $2$-orbifold are
precisely those groups for which at least one of the groups 
$\mathbf{L}$ and $\mathbf{R}$
belong to the set $\{\mathbf{C}_{2n}, \mathbf{D}_{n}\}$.  
In Du\,Val's enumeration, these are the families \#$1$--$19$ 
(including \#$11a$, $1'$, $11'$, $11a'$).   Among fibering subgroups of 
$SO(4)$, families \#$10$, \#$15$, and \#$19$
are maximal; all other groups are subgroups of some member of these families.  
Since we are looking for a lower bound
on the diameters of the orbifolds arising from these families, it suffices to 
examine the maximal families, as follows.

\smallskip

\noindent
$\mathbf{10}$\ $(\mathbf{D}_m/\mathbf{D}_m; \mathbf{D}_n/\mathbf{D}_n)$:  
The orbit of the quaternion 1 in $S^3$ is the union of $\mathbf{C}_{2L}$ and  
the coset  $\mathbf{i} \mathbf{C}_{2L}$,
where $L=\lcm(m,n)$.  The pre-fundamental domain is the same as for 
$(\mathbf{D}_L/\mathbf{D}_L; \mathbf{D}_L/\mathbf{D}_L)$, a 2$L$-prism with 
vertices 
$\frac{1}{\sqrt2}(\cos(\frac{\pi}{2L})  \pm \sin (\frac{\pi}{2L}) \mathbf{k} 
+ \cos (\frac{\pi t}{2L}) \mathbf{i}
 + \sin (\frac{\pi t}{2L}) \mathbf{j} ) $ 
 where $t = 1, 3, ..., 4L-1.$ 
A lower bound for the diameter is 
$\arccos(\frac{\cos(\frac{\pi}{2L})}{\sqrt2})$, which  
is always greater than 
$\arccos(1/\sqrt{2})  = \pi/4$.  The diameter approaches
$\pi/4$ as $L \rightarrow \infty$ and the group approaches the 
corresponding cohomogeneity one action.

\smallskip

\noindent
$\mathbf{15}$\ $(\mathbf{D}_m/\mathbf{D}_m; \mathbf{O}/\mathbf{O}) $:
Let $m \rightarrow \infty$.  The limit group $G$ contains every group in 
this family.
Its identity component $G_0$ can be described as
\[
A(t)=\left( \begin{matrix}
\cos(t) & 0 &0 &-\sin(t)\\
0 & \cos(t)&-\sin(t) & 0\\
0 & \sin(t) &\cos(t) & 0\\
\sin(t)& 0 &0 & \cos(t)\\
\end{matrix} \right), 
\]
where $t \in \rrr$.
In addition, $S^3/G_0 = S^2(\frac{1}{2})$ and $G/G_0 = O^h$,
and so a lower bound on the diameter of $S^3/G$
is $\frac{1}{2} \arccos(\frac{1}{\sqrt{3}})$ (see Table \ref{Ta:2Dfin}).

\smallskip

\noindent
$\mathbf{19}$\ $ (\mathbf{D}_m/\mathbf{D}_m; \mathbf{I}/ \mathbf{I})$:
Similarly, $S^3/G_0 = S^2(\frac{1}{2})$ and $G/G_0 = I^h$,
and so a lower bound on the diameter is 
$\frac{1}{2} \arccos(\frac{\tan(\frac{3 \pi}{10})}{\sqrt3})$.

\subsection{Diameters for the remaining fibering subgroups of $\mathbf{O(4)}$}\label{S:fibO}

The remaining fibering subgroups belong to Du\,Val's families \#$33$--$38$.
Each of these groups is contained in some member of family \#$34$, since
$(\mathbf{L}/\textsc{l};\mathbf{R}/\textsc{r};\phi)$ will be contained
in $(\mathbf{D}_n/\mathbf{D}_n; \mathbf{D}_n/\mathbf{D}_n)^*$ 
if $\mathbf{L}$ and $\mathbf{R}$ are contained in $\mathbf{D}_n$.  
Hence, to find a lower bound on the diameter,  it suffices to examine 
this family.

\smallskip

\noindent
$\mathbf{34}$\ $(\mathbf{D}_n/\mathbf{D}_n; \mathbf{D}_n/\mathbf{D}_n)^*$:  
The orientation-preserving subgroup belongs to family \#$10$ (with $m=n$).  
The number of points in the orbit of the quaternion $1$
under a finite subgroup of $O(4)$ equals the order of the group divided by the 
order of the 
isotropy subgroup at $1$.  
The full group has twice the order of the orientation-preserving subgroup, 
but the order of the isotropy subgroup also doubles, since 
$\mathbf{q} \mapsto \overline{\mathbf{q}}$ 
is an orientation-reversing 
element of any group in family \#$34$, and fixes $1$.  
Hence the orbit of $1$ remains the same after extension, so the 
pre-fundamental domain remains the same.
Consequently, the diameter remains at least 
$\arccos(\frac{\cos(\frac{\pi}{2n})}{\sqrt2} )$; see section \ref{S:fibSO}.

\subsection{Diameters for the nonfibering subgroups of $\mathbf{SO(4)}$}\label{S:nfibSO}

The subgroups $G$ of $SO(4)$ for which the corresponding orbifolds
$S^3/G$ do not admit a fibering over a $2$-orbifold are
precisely those groups for which both groups $\mathbf{L}$ and $\mathbf{R}$
belong to the set $\{\mathbf{T},\mathbf{O},\mathbf{I}\}$.  
In Du\,Val's enumeration, these are groups \#$20$--$32$ 
(including $21'$, $26'$, $26''$, $31'$, $32'$).  
The diameter bound from 1 is sharp when 
it equals $\pi/2$ or $\pi$, since the diameter is greater than $\pi/2$ exactly 
when there is a point fixed by the entire group \cite{B91}, \cite{GM95}.  
Among nonfibering subgroups of $SO(4)$, \#$29$ is
the only group which is maximal with respect to inclusion among finite 
subgroups of $O(4)$; its orbit space is a natural candidate for the minimal 
diameter spherical orbifold.  The other groups are either subgroups of it 
or are contained in a subgroup of $O(4)$ which contains orientation-reversing 
transformations.
We present the groups in order of decreasing
diameter, first considering those whose diameter is 
a rational multiple of $\pi$ and then those whose  
whose diameter is 
an irrational multiple of $\pi$.

\smallskip

\noindent 
$\mathbf{21'}$\ $(\mathbf{T}/\mathbf{C}_1;\mathbf{T}/\mathbf{C}_1)$,
$\mathbf{26'}$\ $(\mathbf{O}/\mathbf{C}_1;\mathbf{O}/\mathbf{C}_1;id)$,
$\mathbf{31'}$\ $(\mathbf{I}/\mathbf{C}_1;\mathbf{I}/\mathbf{C}_1)$:
In all these cases, the quaternion $1$ is fixed by the entire group.  
The diameter is $\pi$.

\smallskip

\noindent
$\mathbf{21}$\ $(\mathbf{T}/\mathbf{C}_2;\mathbf{T}/\mathbf{C}_2)$,
$\mathbf{26}$\ $(\mathbf{O}/\mathbf{C}_2;\mathbf{O}/\mathbf{C}_2)$,
$\mathbf{26''}$\ $(\mathbf{O}/\mathbf{C}_1;\mathbf{O}/\mathbf{C}_1;\xi)$,
$\mathbf{31}$\ $(\mathbf{I}/\mathbf{C}_2;\mathbf{I}/\mathbf{C}_2)$: 
The quaternion $1$ is either mapped to itself or to $-1$ under this group.  
The pre-fundamental domain is bounded by the great sphere perpendicular to 
$1$.  
So the diameter is $\pi/2$.

\smallskip

\noindent 
$\mathbf{22}$\ $(\mathbf{T}/\mathbf{D}_2;\mathbf{T}/\mathbf{D}_2)$,
$\mathbf{27}$\ $(\mathbf{O}/\mathbf{D}_2;\mathbf{O}/\mathbf{D}_2)$: 
The closest images of the quaternion $1$ are 
$\pm\mathbf{i}, \pm\mathbf{j}, \pm\mathbf{k}$; indeed, the 
entire orbit of $1$ is $\mathbf{D}_2$.  
The pre-fundamental domain is a cube with vertices
$(1 \pm\mathbf{i} \pm\mathbf{j} \pm\mathbf{k})/2$.
So a lower bound for the diameter is $\arccos(1/2) = \pi/3$. 

\smallskip

\noindent
$\mathbf{20}$\ $(\mathbf{T}/\mathbf{T};\mathbf{T}/\mathbf{T})$,
$\mathbf{28}$\ $(\mathbf{O}/\mathbf{T};\mathbf{O}/\mathbf{T})$:
The closest images of the quaternion $1$ under these groups are 
$(1\pm \mathbf{i}\pm \mathbf{j} \pm \mathbf{k})/2$ (their 
isotropy subgroups are, respectively, tetrahedral and octahedral). 
The pre-fundamental domain is a
regular octahedron with vertices 
$(1\pm \mathbf{i})/\sqrt{2}, (1\pm \mathbf{j})/\sqrt{2}, 
(1\pm \mathbf{k})/\sqrt{2}$.  
So a lower bound for the diameter is $\arccos(1/\sqrt{2}) = \pi/4$.

\smallskip

\noindent
$\mathbf{32'}$\ $(\mathbf{I}^\dag/\mathbf{C}_1;\mathbf{I}/\mathbf{C}_1;\phi_\dag^{-1})$:  
The closest images of the quaternion $1$ are 
$(-1+\sqrt{5}\mathbf{i}+\sqrt{5}\mathbf{j}+\sqrt{5}\mathbf{k})/4$
plus the 3 additional points obtained by changing two plus signs
to minus signs.  
The pre-fundamental domain is a regular tetrahedron with vertices 
(antipodal to these images)
$(1-\sqrt{5}\mathbf{i}-\sqrt{5}\mathbf{j}-\sqrt{5}\mathbf{k})/4$
plus the 3 additional points obtained by changing two minus signs
to plus signs.  
So a lower bound for the diameter is $\arccos(1/4) \approx \pi/2.38$. 

\smallskip

\noindent 
$\mathbf{32}$\ 
$(\mathbf{I}^\dag/\mathbf{C}_2;\mathbf{I}/\mathbf{C}_2;\tilde{\phi}_\dag^{-1})$:
The closest images of the quaternion $1$ are 
$(1-\sqrt{5}\mathbf{i}-\sqrt{5}\mathbf{j}-\sqrt{5}\mathbf{k})/4$
together with
$(1-\sqrt{5}\mathbf{i}+\sqrt{5}\mathbf{j}+\sqrt{5}\mathbf{k})/4$,
plus the 2 additional points obtained by cyclically permuting 
$\mathbf{i},\mathbf{j},\mathbf{k}$ in the latter expression.  
The next closest images are the 4 points antipodal to these.
The pre-fundamental domain is a truncated regular tetrahedron with vertices 
$(\sqrt{5}+3\mathbf{i}+\mathbf{j}+\mathbf{k})/4$, plus the 11 additional
points obtained by changing two plus signs to minus signs and/or cyclically 
permuting 
$\mathbf{i},\mathbf{j},\mathbf{k}$.  
So a lower bound for the diameter is $\arccos(\sqrt{5}/4) \approx \pi/3.21$.

\smallskip

\noindent 
$\mathbf{23}$\ $(\mathbf{T}/\mathbf{T};\mathbf{O}/\mathbf{O})$,
$\mathbf{25}$\ $(\mathbf{O}/\mathbf{O};\mathbf{O}/\mathbf{O})$:
The closest images of the quaternion $1$ under this group are 
$(1\pm \mathbf{i})/\sqrt{2}, (1\pm \mathbf{j})/\sqrt{2}, 
(1\pm \mathbf{k})/\sqrt{2}$;
next closest are $(1\pm \mathbf{i}\pm \mathbf{j} \pm \mathbf{k})/2$ (from the
subgroup \#$20$).  The pre-fundamental domain is a truncated cube, with vertices 
$((\sqrt{2}+1)\pm(\sqrt{2}-1)\mathbf{i}\pm\mathbf{j}\pm\mathbf{k})/2\sqrt{2}$, 
plus the 16 additional points obtained by cyclically permuting 
$\mathbf{i},\mathbf{j},\mathbf{k}$.
So a lower bound for the diameter is $\arccos((\sqrt{2}+1)/2\sqrt{2}) 
\approx \pi/5.73$. 

\smallskip

\noindent
$\mathbf{24}$\ $(\mathbf{T}/\mathbf{T};\mathbf{I}/\mathbf{I})$,
$\mathbf{30}$\ $(\mathbf{I}/\mathbf{I};\mathbf{I}/\mathbf{I})$:
For group \#$24$, the closest images of the quaternion $1$ under this group are
$ ((\sqrt{5}+1) + 2 \mathbf{i} + (\sqrt{5}-1) \mathbf{j} + 0  \mathbf{k} )/4$, 
plus the eleven additional points  in the orbit of this point under the action 
of the tetrahedral group (the image of $\mathbf{T}$ in  $SO(3)$).
The same is true for group \#$30$, except that 
instead of 11 additional points we have 8 obtained by 
cyclically permuting $\mii,\mj$ and $\mk$.  
The pre-fundamental domain
is a dodecahedron with vertices
$ ( (3\sqrt{2}+\sqrt{10})   + (3+\sqrt{5})(\sqrt{10}-2\sqrt{2})   
\mathbf{i}  +  (3+\sqrt{5})(\sqrt{10}-2\sqrt{2})  \mathbf{j}  + 
(3+\sqrt{5})(\sqrt{10}-2\sqrt{2})  \mathbf{k} ) /8 $,
plus the 7 additional points obtained by reflections in the three 
great spheres orthogonal to 
$\mathbf{i}$, $\mathbf{j}$, and $\mathbf{k}$.  In addition, 
$ ( (6\sqrt{2}+2\sqrt{10}) + 0 \mathbf{i} + (7-3\sqrt{5})(3\sqrt{2}+
\sqrt{10})\mathbf{j}  + 
4\sqrt{2}\mathbf{k}    ) /16 $ is a vertex, as well as the 11 additional 
points in the orbit of this point
under the action of the tetrahedral group (the image of $\mathbf{T}$ in 
$SO(3)$).  So a lower bound for the diameter is
$\arccos((3\sqrt{2}+\sqrt{10} )/8) \approx \pi/8.10$.

\smallskip

\noindent 
$\mathbf{29}$\ $(\mathbf{O}/\mathbf{O}; \mathbf{I}/ \mathbf{I})$: 
The closest images of the quaternion $1$ under this group are
$((3\sqrt{2}+\sqrt{10}) + (\sqrt{10}-\sqrt{2}) \mathbf{i}+ 
(\sqrt{10}-\sqrt{2}) \mathbf{j}+ (\sqrt{10}-\sqrt{2}) \mathbf{k})/8$, 
$((3\sqrt{2}+\sqrt{10}) + (\sqrt{10}-\sqrt{2}) \mathbf{i}- 
(\sqrt{10}-\sqrt{2}) \mathbf{j}- (\sqrt{10}-\sqrt{2}) \mathbf{k})/8$, plus 
the two additional points obtained by cyclically permuting $\mathbf{i}, 
\mathbf{j}$, and 
$\mathbf{k}$ in the latter expression.
The next closest images are $((\sqrt{5}+1)+2\mathbf{i} +(\sqrt{5}-1)\mathbf{j}+
0\mathbf{k})/4$, plus the eleven additional points  in the orbit of this point 
under the action of the tetrahedral group (the image of $\mathbf{T}$ in 
$SO(3)$).  The third layer of images are 
$(\sqrt{10}-\sqrt{2}\mathbf{i}-\sqrt{2}\mathbf{j}-\sqrt{2}\mathbf{k} )/4$ 
together with 
$(\sqrt{10}-\sqrt{2}\mathbf{i}+\sqrt{2}\mathbf{j}+\sqrt{2}\mathbf{k} )/4$, 
plus the two additional points
obtained by cyclically permuting $\mathbf{i}, \mathbf{j}$, and 
$\mathbf{k}$ in the latter expression.

The pre-fundamental domain has 4 twelve-sided faces, 4 six-sided faces, and 
12 faces which are
isosceles triangles.   It can also be described as the intersection of a 
smaller tetrahedron with a 
larger tetrahedron in dual position, all the vertices of which are then 
truncated by the intersection with 
a dodecahedron.  
Refer to Figure 9 in \cite{Du94} for details.  
The vertices of one of the isosceles triangles are as follows:
{\footnotesize
\[
\frac{1 +(3-\sqrt{10})\mathbf{i} +(2+\frac{3}{2}\sqrt{2}-\sqrt{5}-
\frac{1}{2}\sqrt{10})\mathbf{j} + 
(1+\frac{1}{2}\sqrt{2}-\frac{1}{2}\sqrt{10})\mathbf{k}}    
{\sqrt{40+12\sqrt{2}-8\sqrt{5}-12\sqrt{10}}}
\]
\[
\frac{1 +(-1-\frac{1}{2}\sqrt{2}+\frac{1}{2}\sqrt{10})\mathbf{i} +
(4-2\sqrt{2}+\sqrt{5}-\sqrt{10})\mathbf{j} + 
(-5+\frac{7}{2}\sqrt{2}-2\sqrt{5}+\frac{3}{2}\sqrt{10})\mathbf{k}}    
{\sqrt{136-90\sqrt{2}+56\sqrt{5}-42\sqrt{10}}}
\]
\[  
\frac{1 +(-3+3\sqrt{2}-2\sqrt{5}+\sqrt{10})\mathbf{i} +
(-2+\sqrt{2}-\sqrt{5}+\sqrt{10})\mathbf{j} + 
(3-3\sqrt{2}+2\sqrt{5}-\sqrt{10})\mathbf{k}}   
 {\sqrt{136-90\sqrt{2}+56\sqrt{5}-42\sqrt{10}}}
\]}

The first vertex is further from 1 than the other two vertices, which have 
the same distance from 1.  
The remaining 33 vertices are obtained by letting the tetrahedral group act 
on  these three vertices,
which leaves invariant the distances to 1.  Hence a lower bound on the 
diameter is
\[ 
\arccos( 1 /(\sqrt{40+12\sqrt{2}-8\sqrt{5}-12\sqrt{10}}) ) 
\approx \pi/8.93.
\]  

The group is maximal
and the diameter is the smallest achieved by a nonfibering subgroup.

\subsection{Diameters for the remaining nonfibering subgroups of 
$\mathbf{O(4)}$}\label{S:nfibO}

For many of these groups, the orbit of the quaternion $1$ is the same
as its orbit under the orientation-preserving subgroup of index 2.  
This occurs if and only if some orientation-reversing
element fixes $1$, which in turn is equivalent to the condition 
$\mathbf{b} = \mathbf{a}^{-1}$
in the description of the group, and hence is often easy to verify by 
inspection.  In fact, it turns out that whenever an orientation-reversing
element fixes $1$, either there is an element with $\mathbf{a} = 1 = \mathbf{b}$
or there is one with $\mathbf{a} = 1/\sqrt{2} + (1/\sqrt{2})\mathbf{k}$,
$\mathbf{b} = 1/\sqrt{2} - (1/\sqrt{2})\mathbf{k}$.
In these cases, the analysis repeats that of the subgroup, 
so we refer the reader back to that subgroup for information
about the images of $1$ and the pre-fundamental domain. The
only groups we need to consider then are 
\#$39m$, $40m$, $44m$, $46$, $49m$ and $51m$. 

\smallskip

\noindent
$\mathbf{39m}$\ $(\mathbf{T}/\mathbf{C}_1; \mathbf{T}/ \mathbf{C}_1)^*_{c-}$,
$\mathbf{40m}$\ $(\mathbf{T}/\mathbf{C}_1; \mathbf{T}/ \mathbf{C}_1)^*_{-}$,
$\mathbf{44m}$\ $(\mathbf{O}/\mathbf{C}_1; \mathbf{O}/ \mathbf{C}_1;id)^*_{-}$,
$\mathbf{49m}$\ $(\mathbf{I}/\mathbf{C}_1; \mathbf{I}/ \mathbf{C}_1)^*_{-}$:
The quaternion $1$ is mapped to itself by all elements of the orientation-preserving
subgroup 
(\#$21'$ in the first two cases, \#$26'$ in the third case, \#$31'$ in the fourth case) 
and mapped to $-1$ by all other elements.  The 
pre-fundamental domain is bounded by the great sphere perpendicular to $1$,
so the diameter is $\pi/2$ in all four cases.

\smallskip

\noindent 
$\mathbf{46}$\ $(\mathbf{O}/\mathbf{T}; \mathbf{O}/ \mathbf{T})^*_-$: 
The orientation-preserving subgroup is \#28.  
The orbit of the quaternion $1$ is $\mathbf{O}$, and hence 
the pre-fundamental domain is the same truncated cube as for group \#$25$.  
So a lower bound for the diameter is $\arccos((\sqrt{2}+1)/2\sqrt{2}) \approx 
\pi/5.73$. 

\smallskip

\noindent 
$\mathbf{51m}$\ $(\mathbf{I}^\dag/\mathbf{C}_1;\mathbf{I}/\mathbf{C}_1;
\phi_\dag^{-1})^*_{-}$:
The orientation-preserving subgroup is \#$32'$. 
The orbit of the quaternion $1$ is the same as that of \#$32$, 
and hence the pre-fundamental domain is the same truncated tetrahedron.    
So a lower bound for the diameter is $\arccos(\sqrt{5}/4) \approx \pi/3.21$.

\medskip

In addition, we include descriptions of reflection groups and their Coxeter graphs
in Table \ref{Ta:ReflO}, cf.\ \cite {GB71}.  For these groups,
the diameter equals the minimum distance between vertices of the
fundamental domain, a spherical polyhedron, so we can supply the
exact diameter, not just a lower bound. 
The appearance of $\Sigma$ in the column for the Coxeter graph
signifies that the group is a suspension to $O(4)$ of the group
of reflections in $O(3)$ which follows (in other words, the
action is reducible and acts trivially on the extra dimension, as
in the $SO(3)$ action on $S^3$ given near the start of section \ref{S:cohomOne}).
\begin{table}
  \begin{center}
    \caption{Reflection Subgroups of $O(4)$}\label{Ta:ReflO}
    \begin{tabular}{llc}
Du\,Val \# &Coxeter graph & Diameter \\
$40p$ & $\Sigma$
\bpcm(80,40)
\put(5,20){\circle{4}}
\put(25,20){\circle{4}}
\put(45,20){\circle{4}}
\put(7,20){\line(1,0){16}}
\put(27,20){\line(1,0){16}}
\epcm
& $\pi$ \\
$44p$ & $\Sigma$
\bpcm(80,40)
\put(5,20){\circle{4}}
\put(25,20){\circle{4}}
\put(45,20){\circle{4}}
\put(7,20){\line(1,0){16}}
\put(27,20){\line(1,0){16}}
\put(33,25){4}
\epcm
& $\pi$ \\
$49p$ & $\Sigma$
\bpcm(80,40)
\put(5,20){\circle{4}}
\put(25,20){\circle{4}}
\put(45,20){\circle{4}}
\put(7,20){\line(1,0){16}}
\put(27,20){\line(1,0){16}}
\put(33,25){5}
\epcm
& $\pi$ \\
$44mp$ &
\bpcm(80,40)
\put(5,20){\circle{4}}
\put(25,20){\circle{4}}
\put(45,20){\circle{4}}
\put(65,20){\circle{4}}
\put(47,20){\line(1,0){16}}
\put(27,20){\line(1,0){16}}
\epcm
& $\pi/2$ \\
$44$ & 
\bpcm(80,40)
\put(5,20){\circle{4}}
\put(25,20){\circle{4}}
\put(45,20){\circle{4}}
\put(65,20){\circle{4}}
\put(47,20){\line(1,0){16}}
\put(27,20){\line(1,0){16}}
\put(53,25){4}
\epcm
& $\pi/2$ \\
$49$ & 
\bpcm(80,40)
\put(5,20){\circle{4}}
\put(25,20){\circle{4}}
\put(45,20){\circle{4}}
\put(65,20){\circle{4}}
\put(47,20){\line(1,0){16}}
\put(27,20){\line(1,0){16}}
\put(53,25){5}
\epcm
& $\pi/2$ \\
$47$ & 
\bpcm(80,40)
\put(5,20){\circle{4}}
\put(25,20){\circle{4}}
\put(45,20){\circle{4}}
\put(65,20){\circle{4}}
\put(47,20){\line(1,0){16}}
\put(27,20){\line(1,0){16}}
\put(7,20){\line(1,0){16}}
\put(53,25){4}
\epcm
& $\pi/3$ \\
$42$ & 
\bpcm(80,40)
\put(5,30){\circle{4}}
\put(25,30){\circle{4}}
\put(45,30){\circle{4}}
\put(25,10){\circle{4}}
\put(27,30){\line(1,0){16}}
\put(7,30){\line(1,0){16}}
\put(25,12){\line(0,1){16}}
\epcm
& $\pi/3$ \\
$51p$ &
\bpcm(80,40)
\put(5,20){\circle{4}}
\put(25,20){\circle{4}}
\put(45,20){\circle{4}}
\put(65,20){\circle{4}}
\put(47,20){\line(1,0){16}}
\put(27,20){\line(1,0){16}}
\put(7,20){\line(1,0){16}}
\epcm
& $\pi/4$ \\
$45$ &
\bpcm(80,40)
\put(5,20){\circle{4}}
\put(25,20){\circle{4}}
\put(45,20){\circle{4}}
\put(65,20){\circle{4}}
\put(47,20){\line(1,0){16}}
\put(27,20){\line(1,0){16}}
\put(7,20){\line(1,0){16}}
\put(33,25){4}
\epcm
& $\pi/4$ \\
$50$ &
\bpcm(80,40)
\put(5,20){\circle{4}}
\put(25,20){\circle{4}}
\put(45,20){\circle{4}}
\put(65,20){\circle{4}}
\put(47,20){\line(1,0){16}}
\put(27,20){\line(1,0){16}}
\put(7,20){\line(1,0){16}}
\put(53,25){5}
\epcm
&  $\arccos(\frac{3+\sqrt{5}}{4\sqrt{2}})$ \\
    \end{tabular}
  \end{center}
\end{table}
%
\section{Conclusions}\label{S:concl}

We summarize the results for cohomogeneity three in three tables:
Table \ref{Ta:diamFib}, for the fibering groups,  
Table \ref{Ta:diamNonfibRa} for the nonfibering groups with diameters 
a rational multiple of $\pi$, and  
Table \ref{Ta:diamNonfibIr} for the remaining nonfibering groups.
\begin{table}
  \begin{center}
  \caption{Diameters for Finite, Fibering Groups}\label{Ta:diamFib}
    \begin{tabular}{lc}
Du\,Val \#  & Lower Bound for Diameter \\
10 & $\frac{\pi}{4}$ \\
15 & $\frac{1}{2}\arccos(\frac{1}{\sqrt{3}})$\\
19 & $\frac{1}{2}\arccos(\frac{\tan(\frac{3 \pi}{10})}{\sqrt3})$\\
34 & $\frac{\pi}{4}$\\
    \end{tabular}
  \end{center}
\end{table}
\begin{table}
  \begin{center}
  \caption{Diameters in $\pi\mathbb{Q}$ for Finite, Nonfibering Groups}\label{Ta:diamNonfibRa}
    \begin{tabular}{lcc}
Du\,Val \# & Maximal Inclusions & Diameter\\
$21'$ &&  $\pi$\\
$26'$ & $21'\normal 26'$ & $\pi$\\ 
$31'$ && $\pi$\\
$39p$ & $21'\normal 39p$ & $\pi$\\
$40$ & $21\normal 40$ & $\pi$\\
$40p$ & $21'\normal 40p$ & $\pi$\\
$44$ & $26\normal 44$ & $\pi$\\
$44p$ & $26'\normal 44p$ & $\pi$\\
$49p$ & $31'\normal 49p$ & $\pi$\\
\\
$21$ & $21'\normal 21$ & $\frac{\pi}{2}$\\
$26$ & $21\normal 26$ & $\frac{\pi}{2}$\\
 & $26'\normal 26$ &\\
 & $26''\normal 26$ &\\
$26''$ & $21'\normal 26'$ & $\frac{\pi}{2}$\\
$31$ & $31'\normal 31$ & $\frac{\pi}{2}$\\        
$39$ & $21\normal 39$ & $\frac{\pi}{2}$\\
$39m$ & $21'\normal 39m$ & $\frac{\pi}{2}$\\
$40m$ & $21'\normal 40m$ & $\frac{\pi}{2}$\\
$44m$ & $26'\normal 44m$ & $\frac{\pi}{2}$\\
$44pm$ & $26''\normal 44pm$ & $\frac{\pi}{2}$\\
$44mp$ & $26''\normal 44mp$ & $\frac{\pi}{2}$\\
$49$ & $31\normal 49$ &  $\frac{\pi}{2}$\\
$49m$ & $31'\normal 49m$ & $\frac{\pi}{2}$\\
\\
$22$ & $21\subset 22$ & $\frac{\pi}{3}$\\
$27$ & $26\subset 27$ & $\frac{\pi}{3}$\\
 & $22\normal 27$ & \\
$41$ & $22\normal 41$ & $\frac{\pi}{3}$\\
$42$ & $22\normal 42$ & $\frac{\pi}{3}$\\
$47$ & $27\normal 47$ & $\frac{\pi}{3}$\\
 & $26\subset 47$ & \\
\\
$20$ & $22\subset 20$ &  $\frac{\pi}{4}$\\
$28$ & $27\subset 28$ & $\frac{\pi}{4}$\\
 & $20\normal 28$ & \\
$43$ & $20\normal 43$ & $\frac{\pi}{4}$\\
$45$ & $28\normal 45$ & $\frac{\pi}{4}$\\
    \end{tabular}
  \end{center}
\end{table}
\begin{table}
  \begin{center}
  \caption{Diameters Not in $\pi\mathbb{Q}$ for Finite, Nonfibering Groups}\label{Ta:diamNonfibIr}
    \begin{tabular}{lcc}
Du\,Val \# & Maximal Inclusions & Diameter\\
$32'$ && $\arccos(\frac{1}{4})$\\
$51p$ &$32'\normal 51p$ & $\arccos(\frac{1}{4})$\\
&&\\
$32$ &$32'\normal 32$& $\arccos(\frac{\sqrt{5}}{4})$\\
$51$ &$32\normal 51$ & $\arccos(\frac{\sqrt{5}}{4})$\\
$51m$ &$32\normal 51m$ &$\arccos(\frac{\sqrt{5}}{4})$\\
&&\\
$23$ &$20\normal 23$& $\arccos(\frac{\sqrt{2} + 1}{2\sqrt{2}})$\\
$25$ &$23\normal 25$& $\arccos(\frac{\sqrt{2} + 1}{2\sqrt{2}})$\\
& $28\normal 25$&\\
$46$ & $28\normal 46$& $\arccos(\frac{\sqrt{2} + 1}{2\sqrt{2}})$\\
$48$ & $25\normal 48$& $\arccos(\frac{\sqrt{2} + 1}{2\sqrt{2}})$\\
$24$ & $20\subset 24$ & $\arccos(\frac{3\sqrt{2} + \sqrt{10}}{8})$\\
$30$ & $31\subset 30$,&$\arccos(\frac{3\sqrt{2} + \sqrt{10}}{8})$\\
& $g 32 g^{-1}\subset 30$,&\\
& where $g\in SO(4)$ &\\
$50$ & $30\normal 50$ & $\arccos(\frac{3\sqrt{2} + \sqrt{10}}{8})$\\
&&\\
$29$ & $24\normal 29$ & 
{\footnotesize
$\arccos((\sqrt{40} + 12\sqrt{2} - 8\sqrt{5} - 12\sqrt{10})^{-1})$}\\
    \end{tabular}
  \end{center}
\end{table}

All the normal subgroups in these tables have index 2. For the remaining subgroups
we note that the inclusions of \#$22$ in \#$20$ and \#$27$ in \#$28$ are both index 3, 
while the inclusions of \#$21$ in \#$22$ and \#$26$ in \#$27$ are both index 4,
the inclusion of \#$20$ in \#$24$ is index 5, and the inclusions of
\#$31$ and a conjugate in $\so(4)$ of \#$32$ in \#$30$ are both index 60.

We summarize this information in the following theorem:
\begin{theorem}\label{T:MinDiamFin}
Let $G$ be a non-trivial finite subgroup of $O(4)$. Then
\[
\min(\diam(S^n/G)) = \left\{ 
\begin{array}{ll}
\frac{\alpha}{2} & \textrm{for fibering groups}\\
\beta & \textrm{for nonfibering groups}\\
\end{array} 
\right.
\]
where 
$\alpha=\arccos (\frac{\tan(\frac{3 \pi}{10})}{\sqrt3}), 
\beta=\arccos((\sqrt{40} + 12\sqrt{2} - 8\sqrt{5} - 12\sqrt{10})^{-1})$.
\end{theorem}
Note that $\alpha/2 \approx \pi/9.63$ is strictly smaller than $\beta \approx \pi/8.93$.
We further note that there seems to be no relationship whatsoever between
the index of an extension and any subsequent change in diameter. In 
particular, there are many extensions of index 2 where the diameter remains 
unchanged, others where the diameter is reduced by half, others where the 
diameter is reduced instead by $\frac{3}{4}$ and still others by
$\beta/(\arccos(\frac{3\sqrt{2} + \sqrt{10}}{8}))$.

We observe that the nonfibering groups giving ``nice'' diameters, i.e., 
rational multiples of $\pi$, are mainly of one type. The majority 
are of the form $(\mathbf{L}/\textsc{l};\mathbf{R}/\textsc{r})$ 
where $\mathbf{L}=\mathbf{R}\in\{\mt, \mo, \mi\}$
and $\textsc{l}=\textsc{r}\in\{\mc_2,\md_2\}$.
There are two exceptions, groups \#$20$ and \#$28$,  
which are of the form $(\mt/\mt; \mt/\mt)$ and $(\mo/\mt; \mo/\mt)$.
The former groups with $\textsc{l}=\textsc{r}=\mc_2$ are ``diagonal'' groups; 
that is, they project to diagonal 
subgroups of $\so(3)\times \so(3)$. These groups are characterized by
having very small orbits of the quaternion $1$ and hence very large
pre-fundamental domains. If the group contains the nontrivial central element
then the orbit of $1$ is $\pm 1$, and if not, its orbit is simply $1$.
Thus these groups will have pre-fundamental domain the entire 3-sphere
or the half sphere and corresponding diameters of $\pi$ and $\frac{\pi}{2}$.
The bigger $\textsc{l}$ and $\textsc{r}$ become the larger the orbit of $1$ becomes 
(for fixed $\mathbf{L}$, $\mathbf{R}$). These will be points at rational multiples of $\pi$
away from $1$, generating faces of the pre-fundamental domain 
(halfway between each point and $1$). 
However, there is no easy way to predict in general whether or not these faces
will happen to intersect at a point which is at a distance from $1$
which is a rational multiplle of $\pi$.

The remaining nonfibering groups giving diameters that are irrational 
multiples of $\pi$ are generally of the form
$(\mathbf{L}/\textsc{l};\mathbf{R}/\textsc{r})$, 
where $\mathbf{L}, \mathbf{R}\in \{\mt, \mo, \mi\}$ and 
$\textsc{l}, \textsc{r}\in \{\mt, \mo, \mi\}$ or of the 
form ($\mi/\mc_1; \mi/\mc_1)$ or an extension of the same.

Note that $(\mt/\mt; \mt/\mt)\normal (\mt/\mt; \mo/\mo)\normal 
(\mo/\mo; \mo/\mo)$
and $(\mo/\mt; \mo/\mt)\normal (\mo/\mo; \mo/\mo)$, but  
$(\mt/\mt; \mt/\mt)$ and $(\mo/\mt; \mo/\mt)$
both have diameters that are rational multiples of $\pi$
and $(\mt/\mt; \mo/\mo)$, while $(\mo/\mo; \mo/\mo)$ both have diameters that are
irrational multiples of $\pi$. Thus a general theorem relating these finite
groups with diameters that are rational or irrational 
multiples of $\pi$ seems elusive.

\medskip

We are interested in finding a global lower bound for 
isometric group actions on spheres.
We note that the lower bound for any polar action arising from 
a symmetric space $G/H$ where either $G$ or $H$ (or both) 
is a product of classical Lie groups only, which we will define
as a {\it classical polar} action,  approaches $\pi/2$ as 
the cohomogeneity of the action increases (cf. \cite{MS06})). 
Those polar actions arising from symmetric 
spaces for which $G$ or $H$ (or both) is a product of classical Lie groups and
exceptional Lie groups are to be called {\it exceptional polar} actions.
The result holds for many of these groups as well, but for 
many of the groups in this list, the orbit space is yet to be calculated
(given that these admit no ``easy'' matrix representation, other methods 
must be used). We note as well, that in the spherical cohomogeneity 2 case,
the possible orbit spaces are $S^2$, $D^2$ (with 0, 1, 2 or 3 exceptional 
singular points corresponding to an isolated singular orbit).
For the disk cases, those with exceptional singular orbits are 
limited by the corresponding isolated singular orbits as to 
further possible identifications, just as in the case of the interval for 
cohomogeneity 1). In particular, there are very few cases where 
the diameter will actually be changed after an identification and
in none of these cases does the overall minimum diameter
change. For the disk with no exceptional singular points, 
the same holds true.  
Given the work done here and work from \cite{MS05} and 
\cite{MS06} on classical 
connected polar actions of cohomogeneity 3, we can prove the following theorem:
\begin{theorem}\label{T:MinDiamAlln}
Let $G$ act irreducibly 
by cohomogeneity 1, 2, or 3 on $S^n$. Further suppose 
that the action is classical connected polar or non-trivial disconnected 
for the cohomogeneity 3 cases when $n=3$. 
Then
\[
\min(\diam(S^n/G)=\left\{ 
\begin{array}{ll}
\frac{\pi}{6} & \textrm{for cohomogeneity 1}\\
\frac{\alpha}{2} & \textrm{for cohomogeneity 2}\\
\frac{\alpha}{2} & \textrm{for cohomogeneity 3}
\end{array} 
\right.
\]
where 
$\alpha=\arccos (\frac{\tan(\frac{3 \pi}{10})}{\sqrt3})$
\end{theorem}
We further include the following conjecture:
\begin{conjecture}\label{T:limitConj}
Let $G$ act irreducibly on $S^n$ by cohomogeneity $k$, where $n\in 2\zzz$.  Then 
for all $\epsilon > 0$ and for all sufficiently large $k$
(and all $n > k$), $\diam(S^n/G)$ is within $\epsilon$ of $\pi/2$.
\end{conjecture}
We base this conjecture on the following:
in \cite{MS05} and \cite{MS06}, we can show that no further finite 
identifications can be made on the orbit spaces in cohomogeneity 3.
Thus, there are no finite extensions of the classical connected polar actions which will decrease
diameter.  Furthermore, for a cohomogeneity $k$ action on $S^n$,  where $3\leq k<n$ and 
additionally  $k\neq 3, 7$ when $n=7, 15$ respectively,
$S^k$ does not appear as a quotient space, in contrast to the $S^2$ which results from a cohomogeneity 2 circle action on $S^3$. Thus we can reasonably expect diameters to be strictly larger than those found in cohomogeneities 2 and 3,
for sufficiently large $k$. We exclude the odd-dimensional spheres, because
they all admit cohomogeneity $n-1$ free circle actions with quotient spaces a complex projective space
of diameter $\frac{\pi}{2}$ and hence are likely to have quotients by finite extensions with 
diameters approaching $\frac{\pi}{4}$.

When considering cohomogeneity $n$ actions on $S^n$, $n\geq 4$, one might expect these diameters to increase toward $\frac{\pi}{2}$ as well,  at least for $n$ even, since for 
$n$ uneven, we will have once again, as in section \ref{S:fibSO},  finite actions converging to a 
free circle action at least halving the diameter of the corresponding spherical quotient space.
To illustrate this tendency towards $\pi/2$,  consider  
the orbifold
$S^n/G$, where $G$ is the full group of symmetries of the cubic tessellation
of $S^n$ (i.e., the  radial projection of a hypercube in $\mathbb{R}^{n+1}$ onto
$S^n$).  Its diameter is $\arccos(1/\sqrt{n+1})$, which converges to
$\pi/2$ as $n$
goes to infinity.


\end{document}